\documentclass{amsart}

\usepackage{amsbsy,amsfonts,amsmath,amssymb,eucal,mathrsfs,graphicx}
\usepackage[all,cmtip]{xy}
\usepackage[utf8]{inputenc}
\usepackage[T1]{fontenc}
\usepackage{tikz}
\usetikzlibrary{arrows}
\usetikzlibrary{matrix}
\usetikzlibrary{positioning}
\usepackage{hyperref}

\usepackage{setspace} 
\usepackage{ stmaryrd }

\newtheorem{counter}[subsubsection]{$\!\!$}

\newenvironment{demo}{{\flushleft \bf Proof~:}}{\hfill $\square$ \vspace{5mm}}

\newenvironment{rien}{\begin{counter} \rm}{\end{counter}}

\newtheorem{theo}{Theorem}[section]
\newtheorem{prop}{Proposition}[section]
\newtheorem{lemm}{Lemma}[section]
\newtheorem{coro}{Corollary}[section]

\theoremstyle{definition}
\newtheorem*{conv}{Convention}
\newtheorem{defi}{Definition}[section]
\newtheorem{exam}{Example}[section]
\newtheorem{exams}{Examples}[section]
\theoremstyle{remark}
\newtheorem*{rema}{Remark}

\newtheorem{subcounter}[subsection]{$\!\!$}
\newenvironment{defi*}{\begin{subcounter} \rm {\bf Definition.}}{\end{subcounter}}
\newenvironment{defis*}{\begin{subcounter} \rm {\bf Definitions.}}{\end{subcounter}}
\newenvironment{nota*}{\begin{subcounter} \rm {\bf Notation.}}{\end{subcounter}}
\newenvironment{situ*}{\begin{subcounter} \rm {\bf Situation.}}{\end{subcounter}}
\newenvironment{prop*}{\begin{subcounter} {\bf Proposition.}}{\end{subcounter}}
\newenvironment{lemm*}{\begin{subcounter} {\bf Lemma.}}{\end{subcounter}}
\newenvironment{fait*}{\begin{subcounter} {\bf Fait.}}{\end{subcounter}}
\newenvironment{coro*}{\begin{subcounter} {\bf Corollaire.}}{\end{subcounter}}
\newenvironment{conj*}{\begin{subcounter} {\bf Conjecture.}}{\end{subcounter}}
\newenvironment{theo*}{\begin{subcounter} {\bf Theoreme.}}{\end{subcounter}}
\newenvironment{quot*}{\begin{subcounter} {\bf Theoreme \ }}{\end{subcounter}}
\newenvironment{cons*}{\begin{subcounter} \rm {\bf Construction.}}{\end{subcounter}}
\newenvironment{rema*}{\begin{subcounter} \rm {\bf Remarque.}}{\end{subcounter}}
\newenvironment{remas*}{\begin{subcounter} \rm {\bf Remarques.}}{\end{subcounter}}
\newenvironment{exem*}{\begin{subcounter} \rm {\bf Exemple.}}{\end{subcounter}}
\newenvironment{exam*}{\begin{subcounter} \rm {\bf Example.}}{\end{subcounter}}
\newenvironment{exams*}{\begin{subcounter} \rm {\bf Examples.}}{\end{subcounter}}
\newenvironment{cont*}{\begin{subcounter} \rm {\bf Contre-exemples.}}{\end{subcounter}}
\newenvironment{rien*}{\begin{subcounter} \rm}{\end{subcounter}}
\newenvironment{appli*}{\begin{subcounter} \rm {\bf Application.}}{\end{subcounter}}

\def\toto{\rightrightarrows}
\def\too{\longrightarrow}
\def\into{\hookrightarrow}

\DeclareMathOperator{\St}{St}

\DeclareMathOperator{\Hom}{Hom}

\DeclareMathOperator{\GL}{GL}

\DeclareMathOperator{\Sch}{Sch}

\DeclareMathOperator{\id}{id}

\DeclareMathOperator{\Spec}{Spec}

\DeclareMathOperator{\Ens}{Ens}

\DeclareMathOperator{\Supp}{Supp}

\DeclareMathOperator{\R}{R}
\DeclareMathOperator{\D}{D}
\DeclareMathOperator{\Gl}{GL}
\DeclareMathOperator{\M}{M}
\DeclareMathOperator{\G}{G}

\DeclareMathOperator{\fix}{fix}

\DeclareMathOperator{\F}{F}

\def\cA{{\mathcal A}} \def\cB{{\mathcal B}} 
 \def\cE{{\mathcal E}} \def\cF{{\mathcal F}}
\def\cG{{\mathcal G}} \def\cH{{\mathcal H}} \def\cI{{\mathcal I}}
  
  \def\cO{{\mathcal O}}
 \def\cQ{{\mathcal Q}} 
 \def\cT{{\mathcal T}}

\renewcommand\AA{\mathbb{A}}

\newcommand\GG{\mathbb{G}}

 \newcommand\ZZ{\mathbb{Z}}

\def\fm{{\mathfrak m}}

\def\fD{{\mathfrak D}}

\def\fS{{\mathfrak S}}

\begin{document}
 
\title[Ramification of inseparable coverings]{Ramification of inseparable coverings of schemes and application to diagonalizable group actions.}
\date{}
\author{Gabriel Zalamansky}
\address{Mathematisch Instituut, Universiteit Leiden, Niels Bohrweg 1,
2333 CA Leiden, The Netherlands.}
\email{g.s.zalamansky@umail.leidenuniv.nl}


\begin{abstract}
 We define the notion of inseparable coverings of schemes and we propose a ramification formalism for them, along the lines of the classical one. Using this formalism we prove a formula analogous to the classical Riemann-Hurwitz formula for generic torsors under infinitesimal diagonalizable group schemes.
\end{abstract}
 
\maketitle


\section{Introduction}

When $f : Y \too X$ is a ramified cover of a smooth scheme $X$, ie a finite, surjective, locally free morphism of smooth schemes which is \'etale over a dense open subscheme of $X$, the classical ramification theory associates to $f$ a divisor that measures the obstruction for $f$ to be an \'etale covering. Let us briefly recall it.
\label{classic}
If $f$ is as above, the sheaf of first-order differential forms $\Omega^1_f$ is trivial on a dense open subscheme of $Y$ and hence is a torsion sheaf, to which one can associate a divisor $\R_f$, also denoted $\R_{Y/X}$, (by a process we recall in \ref{div}) which measures the obstruction for $f$ to be \'etale everywhere on~$Y$. Such a morphism is classically called a ramified covering and the divisor $\R_f$ is called the ramification divisor of $f$. A crucial feature of this construction is that it is transitive with respect to d\'evissage : if $f : Z \too X$ is a ramified covering that factors into ramified coverings $g : Z \too Y$ followed by $h : Y \too X$ then, as divisors on $Z$, we have \[ \label{trans} \tag{*} \R_{Z/X}= \R_{Z/Y} + g^*\R_{Y/X}. \] 
The ramification theory of local rings with perfect residue fields allows for the computation of the local multiplicities of the divisors. For a ramified cover $f : Y \too X$, one can relate $\R_{f}$ to the geometry of the morphism $f$ via the formula \[ \label{RH} \tag{RH} \det(\Omega^1_f) = \cO_Y(\R_f) \] from which is derived, in the case of projective curves, the famous Riemann-Hurwitz formula.
  
Observe in particular that if $f$ is the quotient morphism of $Y$ by the action of a finite \'etale group scheme $G$, then $f$ is \'etale if and only if the action of $G$ is free everywhere on~$Y$, if and only if $f : Y \too X$ is a $G$-torsor. In this case the ramification divisor $\R_f$ also measures the obstruction for $f$ to be a $G$-torsor and the formula \eqref{RH} relates the action on~$Y$ to the geometry of the quotient morphism. 

We now raise the question : what if $G$ is no longer assumed to be \'etale ? 
More precisely, if $Y$ is a scheme of characteristic $p > 0$ and $G$ is a finite flat group scheme (possibly infinitesimal) acting on $Y$, freely on a dense open subset, can one measure the obstruction for the quotient morphism to be a $G$-torsor ? Our goal is thus to develop a theory of ramified coverings in which the unramified objects would no longer be the \'etale morphisms but the torsors. Note that torsors under infinitesimal group schemes are purely inseparable. In this case the sheaf of differential $1$-forms is no longer torsion and one cannot hope to directly carry over the previous definitions of ramification to this setting. We then have to find a substitute for sheaf of differential 1-forms.

There is, however, an issue which lies in the very formulation of these questions that needs to be addressed first, as illustrated by the following example.

\begin{exam*}
\label{ex}
Let $k$ be a field of characteristic $p>0$ and $\AA^1 = \Spec(k[x])$ be the affine line over $k$. Consider the action of the infinitesimal $k$-group scheme $\mu_{p,k} = \Spec(\frac{k[s]}{s^p-1})$ on $\AA^1$ given by $s.x= sx$. An easy computation of the invariant ring shows that the quotient morphism is the absolute Frobenius : \[ \begin{tabular}{cccc}
$\F$ : & $\AA^1$ & $\too$ & $\AA^1$ \\ 

& $x$ & $\mapsto$ & $x^p$ \\ 

\end{tabular}. \] 

This action is free on the dense open subscheme $\AA^1 \setminus \{0\} $ and has a fixed point in $0$. The quotient morphism is thus not a $\mu_{p,k}$-torsor.

However, it is easily seen that $\F$ is also the quotient morphism for the action on $\AA^1$ of the infinitesimal $k$-group scheme $\alpha_{p,k} = \Spec(\frac{k[t]}{t^p})$ given by $t.x = x + t$, which is free everywhere. This makes $\F$ into an $\alpha_{p,k}$-torsor.

Finally, $\F$ can also be seen as the quotient morphism for the non-free action of $\alpha_{p,k}$ on~$\AA^1$, this time given by $t.x = \frac{x}{1+tx} = 1 -tx +...+ (-1)^{p-1}t^{p-1}x^{p-1}$. Note that the group schemes $\alpha_{p,k}$ and $\mu_{p,k}$ are not isomorphic.

\end{exam*}
This example shows that, contrary to classical ramified coverings, neither the group acting nor the eventuality of being a torsor is determined by the sole quotient morphism. In this situation, the question of measuring the obstruction of a finite locally free morphism to be a torsor makes sense only relatively to a given group action.  

\bigskip
\noindent{\textbf{Overview of the paper}}

\medskip

In the first section of this article, we define the notion of "generalized coverings" which includes the data of a specific action along with a finite flat morphism. Our definitions are formulated in terms of groupoid schemes. This allows our formalism to include finite flat morphisms arising from quotients by vector fields (ie. foliations) that do not necessarily stem from group actions. We recall the relevant basic facts about groupoid schemes.
 
In a second section we proceed to propose a ramification formalism along the lines of the classical one that we outlined. In the case of a generically \'etale Galois covering, we recover the ramification divisor of the classical theory.

In the last section we then specify the situation to actions of infinitesimal diagonalizable group schemes to obtain a formula relating the action and the geometry of the quotient morphism, much like the classical Riemann-Hurwitz formula.

\newpage

\section{Groupoids and generalized covers.}

Let us recall the basic facts and properties of groupoid schemes that will be needed to formulate our ramification formalism. Let us fix a base scheme $S$.

\subsection{Groupoid schemes.}

\begin{defi}

A groupoid scheme over $S$ is a quintuplet $(X,\cG,s,t,c)$ where 

\begin{itemize}

\item $X$ and $\cG$ are $S$-schemes,

\item $s$ and $t$ are morphisms of $S$-schemes $\cG \too X$ respectively called source and target,

\item $c : \cG \times_{s,t} \cG \too \cG$ is an $S$-morphism called composition,

\end{itemize}

such that for all $S$-scheme $T$, the $T$-points $(X(T),\cG(T),s,t,c)$ form a category whose objects are elements of $X(T)$ and whose arrows are elements of $\cG(T)$ in which every arrow is invertible. 
There is an unit section $X \too \cG$ mapping an object $x$ to the arrow $\id_x$ and an inverse morphism $\cG \too \cG$ mapping each arrow to its inverse. Both are determined by $s,t$ and $c$.

A groupoid is said to be finite locally free of rank $n$ if the morphism $s$ (or equivalently $t$) is. In that case we use the notation $[\cG : Y] = n$.

A morphism of groupoids $(X,\cG,s,t,c) \too (X',\cG',s',t',c')$ is a map $\cG \too \cG'$ that induces, for all $S$-scheme $T$ , a functor between the two categories of $T$-points. Equivalently, it is a morphism of schemes $f : \cG \too \cG'$ such that the relevant diagram involving $c$ and $c'$ commute.
Note that $f$ induces a morphism on object $f_0 : X \too X'$ defined by $f_0(x) = s'(f(\id_x))$.

 We will often use the notation $\cG \toto X$ for a groupoid $(X,\cG,s,t,c)$. We shall also denote multiplication and inverse multiplicatively.
\end{defi}

\begin{exam}

If $G$ is an $S$-group scheme acting on an $S$-scheme $X$ via a map $a : G \times_S X \too X$ then we get an groupoid scheme $G \times_S X \toto X$ with source $\textup{pr}_2$ and target $a$, often called the action groupoid of $G $ on $X$.

\end{exam}

We extend this terminology to arbitrary groupoids : if $\cG \toto X$ is an $S$-groupoid, we shall say that $\cG$ acts on $X$. The action is said to be \textit{free} if the morphism $j = (s,t) : \cG \too X \times_S X$ is a closed immersion.  

\begin{rien}{\bf{Subgroupoids.}}\end{rien}

If $\cG \toto X$ is an $S$-groupoid, a subgroupoid (resp. closed subgroupoid, resp. open subgroupoid) of $\cG$ is an $S$-groupoid $\cH \toto X$ with an immersion (resp. closed immersion, resp. open immersion) that is a morphism of groupoids. This means that the groupoid structure on $\cH$ is induced by that of $\cG$.

\begin{rien}{\bf{Products.}}\end{rien}

If $(X, \cG,s,t,c)$ and $(X',\cG',s',t',c')$ are two groupoids one can define their product $(X,\cG'',s'',t'',c'')$ in the category of groupoids of objects $X$ as follows : 

\begin{itemize}

\item Set $\cG'' = \cG \times_{X \times_S X} \cG$.

\item Source and target are respectively given by $\textup{pr}_1 \circ s$ and $\textup{pr}_2 \circ t$.

\item Composition is done component-wise, whenever it makes sense.
\end{itemize}
Both projections $\cG'' \too \cG$ and $\cG'' \too \cG'$ are morphisms of groupoids.

\begin{rien}{\bf{Pullbacks.}}\end{rien}

Let $f : Y \too X$ be a morphism between $S$-schemes. If $\cG \toto X$ is a groupoid one can define its pullback to $Y$ as follows : 

\begin{itemize}

\item Set $f^*\cG = (Y\times_S Y) \times_{X \times_S X} \cG$.
\item Set $f^*s(y_1,y_2,g) = y_1$ and $f^*t(y_1,y_2,g) = y_2$.

\item Composition is given by $f^*c((y_1,y_2,g),(z_1,z_2,h)) = (y_1,z_2,gh)$.

\end{itemize}

One can check that $(Y,f^*\cG,f^*s,f^*t,f^*c)$ is an $S$-groupoid which we call the pullback of $\cG \toto X$ by $f$. We shall often denote it by $\cG_{|Y}$.

\begin{rien}{\bf{Kernels.}}\end{rien}

Let $f$ be a morphism between two $S$-groupoids $\cG \toto X$ and $\cG' \toto X'$. The kernel of $f$, denoted by $\ker(f)$ is the groupoid defined as follows : 

\begin{itemize}

\item The schemes of arrows is defined by the fibre product \[ \xymatrix{ \ker(f) \ar[r] \ar[d] & X' \ar[d]^{e'} \\ \cG \ar[r]^f & \cG'} \]
where $e'$ is the unit section of $\cG'$.
\item Source and target are given by the compositions $\ker(f) \too \cG \toto X$.

\item Composition is given by the composition in $\cG$.

\end{itemize} 

By definition points of $\ker(f)$ are those of $\cG$ which are sent to identities by $f$. Let us note that the unit section is an immersion (since $s' \circ e' = \id_{X'}$). Hence $\ker(f)$ is a subgroupoid of $\cG \toto X$.

\begin{rien}{\bf{Stabilizers.}}\end{rien}

Let $\cG \toto X$ be an $S$-groupoid. We define its stabilizer, which we denote $\St_{\cG}$, by the fibre product of $j = (s,t) : \cG \too X \times_S X$ with the diagonal morphism of $X$ : 

\[ \xymatrix{ \St_\cG \ar[r] \ar[d] & \cG \ar[d]^j \\ X \ar[r]^-{\Delta_X} & X \times_S X } \]

The points of $\St_\cG$ are the points of $\cG$ whose source and target are equal. The composition in $\cG$ induces a morphism $\St_\cG \times_X \St_\cG \too \St_\cG$ which makes $\St_\cG$ into an $X$-group scheme. It is the biggest subgroupoid of $\cG$ that is an $X$-group scheme. 
The action of a groupoid is said to be free if it has trivial stabilizer.

\subsection{Generalized covers.}

We now proceed to give a definition of coverings that would include inseparable morphisms invariant under the generically free action of a finite locally free groupoid scheme.

As explained in the introduction, the same inseparable morphism can be seen as the quotient morphism for several actions of non-isomorphic group schemes. Hence to treat these morphisms as coverings we need to specify a groupoid acting on the source.

We propose the following :

\begin{defi}
\label{def}

Fix a base scheme $S$ and an $S$-scheme $X$. A (generalized) covering of $X$ is a couple $(Y \too X, \cG \toto Y)$, where 

\begin{itemize}

\item $Y$ is an $S$-scheme.
\item $\cG \toto Y$ is a finite locally free $X$-groupoid whose orbits are included into open affines of $Y$ and whose action on $Y$ is generically free. This means that there exists a dense open subscheme $V \subset Y$ such that $\cG_{\vert V}$ acts freely on $V$.

\item $Y \too X$ is a finite surjective locally free morphism which is $\cG$-invariant. 

\item The order of $\cG$ is the same as the order of $Y \too X$, ie $[\cG : Y] = [Y:X]$.

\end{itemize}

The word "generalized" will be mostly be employed to stress the difference between classical generically \'etale morphisms and the objects defined above. When no confusion is likely, the latter will just be called coverings.
For short, we shall often write $(Y,\cG)$ instead of $(Y \too X, \cG \toto Y)$. In case $\cG = G \times_S Y$ is the action groupoid for the action of a group $G$, we shall even denote the covering $(Y,G)$. We refer to those as $G$-coverings. 
\end{defi}

Let us note that, if $(Y,\cG)$ is a covering of a scheme $X$, the hypothesis that the orbits of $\cG \toto Y$ are included into open affines of $Y$ imply, together with local freeness, that the quotient $Y/\cG$ exists in the category of $S$-schemes. See \cite[Exp V, th.4.1]{SGA3} for a proof. However, since the action of $\cG$ is not free, the quotient scheme does not represent the fppf quotient sheaf of $Y$ by $\cG$. Hence a priori we do not know the points of $Y/\cG$.

The conditions that we imposed in the definition of a covering imply that the quotient $Y/\cG$ identifies with $X$, according to the following lemma.
\begin{lemm}

Let $(Y,\cG)$ be a covering of an $S$-scheme $X$. The categorical quotient $Y/\cG$ identifies with $X$.

\end{lemm}

\begin{demo}

This follows from the fact that the morphism \[ j_X=(s,t) : \cG \too Y \times_X Y \] is an epimorphism of schemes. Assume this for the moment. We will show that $X$ satisfies the universal property of the quotient $Y/\cG$.

Let $f : Y \too T$ be a $\cG$-invariant morphism of $S$-schemes. Since $Y \too X$ is faithfully flat, to show that $f$ factors through $X$ it suffices, by descent, to show that $f \circ \textup{pr}_1 = f \circ \textup{pr}_2$, where $\textup{pr}_1, \textup{pr}_2 : Y \times_X Y \too Y$ are the two projections. Since by assumption $j_X$ is an epimorphism, it is equivalent to show that $  f \circ \textup{pr}_1 \circ j_X= f \circ \textup{pr}_2 \circ j_X$. But this last equality is just the equality $f \circ s = f \circ t$, which is verified since $f$ is $\cG$-invariant.

Hence we are left to show that $j_X$ is an epimorphism. Note that it is finite, so in particular quasi-compact and quasi-separated. Hence $j_X$ is schematically dominant if and only if $j_X^\sharp : \cO_{Y \times_X Y} \too j_{X*}\cO_\cG$ is injective. Note also that a finite schematically dominant morphism is surjective by Cohen-Seidenberg's theorem. Hence, by \cite[Exp.VIII, Prop 5.1]{SGA1}, to show that $j_X$ is an epimorphism we only have to show that it is schematically dominant.

Let $V \subset Y$ be a saturated schematically dense open subscheme of $Y$ on which $\cG$ acts freely. Since $Y \too X$ is faithfully flat, its image $W \subset X$ is a schematically dense open of $X$. 
By \cite[Exp V, th.4.1]{SGA3}, the morphism $V \too V/\cG$ is finite flat of degree $[Y:X]$. Since $Y \too X$ is $\cG$-invariant, we have a commutative diagram \[ \xymatrix{ V \ar[d] \ar[r] & W \\ V/\cG \ar[ur] } \] from which we conclude, by the fiberwise criterion for flatness, that $V/\cG \too W$ is finite flat of degree $1$, hence an isomorphism. 

Since the action of $\cG_{\vert V}$ is free on $V$, the morphism $j_V : \cG_{\vert V} \too V \times_W V$ is an isomorphism. By faithfull flatness of $Y \too X$, the immersion $V \times_W V \too Y \times_X Y$ is schematically dominant and we have a commutative diagram \[ \xymatrix{\cG_{\vert V} \ar[r] \ar[d]_{j_V} & \cG \ar[d]_{j_X} \\ V \times_W V \ar[r] & Y \times_X Y} \]

from which we conclude that $j_X$ is schematically dominant.

\end{demo}

We can then define the notion of morphisms of generalized coverings, in an obvious way.

\begin{defi}

If $(Y_1,\cG_1)$ and $(Y_2,\cG_2)$ are two coverings of an $S$-scheme $X$, a morphism of coverings is a groupoid morphism $f : \cG_1 \too \cG_2$ such that the following diagram commutes : 

 \[ \xymatrix{ \cG_1 \ar[rr]^{f} \ar@<-2pt>_{s_1}[d] \ar@<2pt>^{t_1}[d] & & \cG_2 \ar@<-2pt>_{s_2}[d] \ar@<2pt>^{t_2}[d] \\  Y_1 \ar[rr]^{f_0} \ar[rd] & & Y_2 \ar[ld] \\ &  X  } \]

\end{defi}

We wish to develop a ramification theory for these objects, along the lines of the classical one, in which the unramified objects would be the coverings given by groupoids acting freely. 

\begin{defi}

A covering $(Y,\cG)$ of an $S$-scheme $X$ is said to be unramified if the groupoid $\cG \toto Y$ acts freely. 

\end{defi}

\section{A ramification divisor for generalized coverings.}

If $(Y,\cG)$ is a covering of an $S$-scheme $X$, by definition the action of $\cG \toto Y$ is free on a dense open subscheme of $Y$. Thus its stabilizer group scheme $\sigma : \St_\cG \too Y$ is trivial over a dense open subscheme of $Y$. Denote by $\fm_\cG$ its augmentation ideal, ie the ideal defining the unit section $Y \too \St_\cG$. The latter is zero if and only if $\St_\cG$ is trivial. Hence it follows that the sheaf of $\cO_Y$-modules $\sigma_*\fm_\cG$ is a torsion sheaf which is trivial if and only if the groupoid $\cG \toto Y$ acts freely.

We thus have a torsion sheaf of $\cO_Y$-modules which is zero exactly when the covering $(Y,\cG)$ is unramified. Accordingly, it is a natural candidate to replace the sheaf of differential $1$-forms of the classical theory.  

We wish to have a geometric incarnation of this sheaf. We use a construction of Mumford, which we recall, that produces an effective Cartier divisor out of a torsion sheaf. Over a smooth scheme (or at least regular in codimension $1$) the corresponding Weil divisor is just the sum of the codimension $1$ points of its support with appropriate multiplicities.

\subsection{Div of a coherent torsion sheaf}
\label{div}
In this section we recall the construction of Mumford that associates an effective Cartier divisor to a coherent torsion sheaf. We refer to \cite[Chap. V.3]{GIT} for greater details.

Let $X$ be a noetherian $S$-scheme and $\cF$ a coherent sheaf on $X$ such that :

\begin{itemize}

\item[(i)] The support of $\cF$ does not contain any associated point (ie depth $0$ point) of $X$.

\item[(ii)] For all point $x \in X$, the stalk $\cF_x$ is of finite tor-dimension, ie admits a finite projective resolution.

\end{itemize}

If $\cE$ is a locally free sheaf of rank $r$ on $X$, we denote by $\det(\cE)$ the invertible sheaf~$\Lambda^r \cE$. Let us start with the following lemma, proved in \cite[Chap.V, \S 3, Lemma 5.6]{GIT}.

\begin{lemm}
If $0 \too \cE_n \too \cE_{n-1} \too ... \too \cE_0 \too 0$ is an exact sequence of locally free sheaves on $X$, there exists a canonical isomorphism $\displaystyle\bigotimes\limits_{i=0}^n \det(\cE_i)^{(-1)^i} \simeq \cO_X$.

\end{lemm}

By assumptions, every point $x \in X$ has an open neighborhood $U$ over which $\cF$ has a finite resolution by free $\cO_U$-modules

\[0 \too \cE_n \too ... \too \cE_0 \too \cF_{\vert U} \too 0. \]

Set $U' = U \setminus \Supp(\cF)$. By definition over $U'$ the sequence \[ 0 \too \cE_{n_{\vert U'}} \too ... \too \cE_{0_{\vert U'}} \too 0 \] is exact, hence by the above lemma there is a canonical isomorphism $\cO_{U'} \simeq \displaystyle\bigotimes\limits_{i=0}^n \det(\cE_i)_{|U'}^{(-1)^i}$.

Also, since the sheaves $\cE_i$ are free on $U$, we have an isomorphism \[\displaystyle\bigotimes\limits_{i=1}^n \det(\cE_i)^{(-1)^i} \simeq \cO_U,\] unique up to a unit. Composing these we get a morphism $\cO_{U'} \too \cO_{U'}$, defined by a section $f \in \cO_X(U')$. Since $f$ is unique up to unit in $U$ and not a zero-divisor by assumption $(i)$, we get a Cartier divisor $(f)$ in $U$. We refer to \cite[Chap. V.3]{GIT} for a proof that these constructions glue to give an effective Cartier divisor $\textup{div}(\cF)$ on~$X$.

\medskip

If $x \in X$ is a point of depth $1$, it follows from the Auslander-Buschbaum formula that, over some neighborhood of $x$,  such a sheaf $\cF$ has a free resolution of the form $0 \too \cE_1 \too \cE_0 \too \cF \too 0$, where $\cE_1$ and $\cE_0$ are free sheaves of the same rank. If $h$ denotes the map $\cE_1 \too \cE_0$, we then have $\textup{div}(\cF)_x = (\det(h))_x$.

When $X$ is regular in codimension $1$, this allows one to give a simple expression of $\textup{div}(\cF)$ :

\begin{lemm}

Suppose $X = \Spec(A)$ is the spectrum of a discrete valuation ring. Let $\pi \in A$ be a uniformizer. There exists an $A$-module $M$ of finite length such that $\cF = \tilde{M}$ and we have $\textup{div}(\cF) = (\pi^{l_A(M)})$, where $l_A(M)$ is the $A$-length of $M$. 

\end{lemm}

\begin{demo}

There exists an isomorphism of $A$-modules $M \simeq \displaystyle\bigoplus\limits_{i=1}^r A/ \pi^{n_i}$ for an $r$-uple of integers $(n_1,...,n_r)$. We then have a resolution \[ 0 \too A^r \stackrel{h}{\too} A^r \too M \too 0, \] 
where $h$ is the diagonal matrix $(\pi^{n_i} \delta_{ij})_{1 \leq i,j \leq n}$ whose determinant is $\pi^{l_A(M)}$.

\end{demo}

This results globalizes immediately to any scheme that is regular in codimension $1$, for if $x \in X$ has codimension $1$ and $U = \Spec(A)$ is an affine neighborhood of $x \in X$, the local ring $\cO_{X,x}$ is a flat $A$-module and we can tensor the resolutions used to compute $\textup{div}(\cF)$ by $\cO_{X,x}$ to obtain resolutions that compute $\textup{div}(\tilde{\cF_x})$. Hence the multiplicity of $\textup{div}(\cF)$ at $x$ is $l_{\cO_{X,x}}(\cF_x)$.

\subsection{A ramification divisor for generalized coverings}

We are now ready to define a ramification divisor for generalized coverings.

\begin{defi}
\label{newdiv}

Let $(Y,\cG)$ be a covering of an $S$-scheme $X$. Let $\sigma : \St_\cG \too Y$ be the stabilizer group scheme of the groupoid $\cG$ and $\fm_\cG$ be its augmentation ideal. Suppose that the $\cO_Y$-module $\sigma_* \fm_\cG$ has finite projective dimension. Define the ramification divisor of $(Y,\cG)$ to be \[ \R_\cG := \textup{div}(\sigma_* \fm_\cG). \]

\end{defi}

\begin{rema}

In order the use the general construction of \ref{div} we have to make sure that the $\cO_Y$-module $\sigma_* \fm_\cG$ has finite projective dimension. This will always be the case if $Y$ is a regular scheme. 

\end{rema}

Let us compute the ramification divisors of the examples \ref{ex} given in the introduction.

\begin{exams}

Let $k$ be a field of caracteristic $p > 0$.

\begin{itemize}

\item[•] Consider the groupoid $\cG \toto \AA^1_k$ given by the action of the group scheme $G = \mu_{p,k} = \Spec( \frac{k[s]}{s^p -1} )$ on the affine line $\AA^1_k = \Spec(k[y])$ by multiplication. It is defined by the coaction 

\[ \begin{tabular}{ccc}
 
$k[y]$ & $\too$ & $\frac{k[y,s]}{s^p-1}$ \\ 
 
$y$ & $\mapsto$ & $sy$ \\ 
 
\end{tabular}. \]

The stabilizer group scheme is given by the fiber product \[ \xymatrix{\St_\cG \ar[r] \ar[d] & G \times_k Y \ar[d]^-j \\ Y \ar[r]^-{\Delta_Y} & Y \times_k Y.} \] 

The morphism $j$ is defined by the ring map \[ \begin{tabular}{ccc}

$k[y_1,y_2]$ & $\too$ & $\frac{k[y,s]}{s^p-1}$ \\ 

$y_1$ & $\mapsto$ & $y$ \\ 

$y_2$ & $\mapsto$ & $sy$ .\\
\end{tabular} \]

Thus we have $\cO_{\St_\cG} = \frac{k[y,s]}{s^p-1} \otimes_{k[y_1,y_2]} k[y] = \frac{k[y,s]}{s^p-1,(s-1)y}$. Its augmentation ideal is $\fm_\cG = (s-1)\frac{k[y,s]}{s^p-1,(s-1)y}$, for which we have the following free resolution as a $k[y]$-module :

\[ 0 \too k[y]^{\oplus p-1} \stackrel{\times y}{\too}k[y]^{\oplus p-1} \too \fm_\cG \too 0 ,\] the first arrow being the multiplication of each coordinates by $y$. Its determinant is $y^{p-1}$. Hence the ramification divisor of this covering is supported in $0 \in \AA^1_k$ where is has multiplicity $p-1$. As Weil divisors we thus have \[ \R_\cG = (p-1) [0]. \]

\item[•] Consider this time the groupoid $\cG \toto \AA^1_k$ given by the action of $G = \alpha_{p,k}$ on the affine line defined by the algebra map
\[ \begin{tabular}{ccc}
 
$k[y]$ & $\too$ & $\frac{k[y,t]}{t^p}$ \\ 
 
$y$ & $\mapsto$ & $\frac{y}{1+ty} = 1 -ty +...+ (-1)^{p-1}t^{p-1}y^{p-1}$ \\ 
 
\end{tabular}. \]

The morphism $j$ is defined by the ring map \[ \begin{tabular}{ccc}

$k[y_1,y_2]$ & $\too$ & $\frac{k[y,t]}{t^p}$ \\ 

$y_1$ & $\mapsto$ & $y$ \\ 

$y_2$ & $\mapsto$ & $\frac{y}{1+ty}$ .\\
\end{tabular} \]

Observe that $y-\frac{y}{1+ty} = \frac{ty^2}{1+ty}$. Hence $\cO_{\St_\cG} = \frac{k[y,t]}{t^p,ty^2}$ and $\fm_\cG = t\frac{k[y,t]}{t^p,ty^2}$.
We have the following free resolution of $\fm_\cG$ as a $k[y]$-module :

\[ 0 \too k[y]^{\oplus p-1} \stackrel{\times y^2}{\too}k[y]^{\oplus p-1} \too \fm_\cG \too 0 ,\] the first arrow being the multiplication of each coordinates by $y^2$. Its determinant is $y^{2(p-1)}$. We thus have \[\R_\cG = 2(p-1)[0].\]

\item[•] Finally, it clear that the action of $\alpha_{p,k}$ by translation on the affine line leads to a groupoid with trivial stabilizer, hence no ramification divisor.
\end{itemize}

We see that our definition of covering allows one to differentiate between these group actions, which was impossible with the sole quotient morphism.
\end{exams}

In case the covering is given by the action of a finite \'etale group scheme, the classical theory already produces a ramification divisor, using first-order differential forms, as recalled in the introduction \ref{classic}. In the next section we show that, in that case, the latter agrees with the one we just defined.

\subsection{The case of generically \'etale Galois coverings.}

\begin{theo}

Let $f : Y \too X$ be a generically \'etale morphism of normal schemes. Suppose that $f$ is a Galois cover of group $G$, in the sense of \cite{SGA1}, and that all the residue fields extensions $k(y)/k(f(y))$ are separable. Denote by $\cG \toto Y$ the action groupoid of $G$ on $Y$. Let $\R_\cG$ be the divisor associated with the stabilizer of $\cG$, defined in \ref{newdiv}, and $\R_{Y/X} = \textup{div}(\Omega^1_{Y/X})$ be the ramification divisor of the classical theory.

One has the equality \[ \R_\cG = \R_{Y/X}. \]

\end{theo}

\begin{demo}

We will show that both divisors have the same multiplicity in each codimension $1$ point of $Y$. By assumption $Y$ is regular in codimension $1$, so we may assume that $Y=\Spec(A)$ and $X=\Spec(A_0)$ are discrete valuation rings, whose corresponding extension $K/K_0$ of fraction fields is Galois of group $G$.

If $B$ is an $A_0$-algebra we shall denote by $B^G$ the algebra of functions from $G$ to $B$, which is a finite $B$-module, a basis being given by the functions  \[ \begin{tabular}{cccc}

$e_g$ :  & $G$ & $\too$ & $B$ \\ 

 & $h$ & $\mapsto$ & $\delta_{h,g}$ \\ 

\end{tabular}. \]

The multiplication in $B^G$ is given by $e_g e_{g'} = \delta_{g,g'}e_g$.

The action of $G$ on $Y$ is given by algebra automorphisms $g^\sharp : A \too A$, one for each $g \in G$, satisfying the usual conditions. If we abuse notations and denote by $g$ the automorphism $(g^\sharp)^{-1}$, the action map $\rho : G \times_X Y \too Y$ corresponds to the algebra map 
  \[ \begin{tabular}{cccc}
 
$\rho^\sharp$ :  & $A$ & $\too$ & $A \otimes_{A_0} A_0[G] \simeq A[G]$ \\ 
 
 & $a$ & $\mapsto$ & $\displaystyle\sum_{g \in G} g(a)e_g$ \\ 
 
\end{tabular}. \]

The morphism $j : G \times_X Y \too Y \times_X Y$ is then given by 

\[ \begin{tabular}{cccc}

$j^\sharp$ :  & $A \otimes_{A_0} A$ & $\too$ & $A[G]$ \\ 

 & $a \otimes b $ & $\mapsto$ & $\displaystyle\sum_{g \in G} ag(b) e_g$ \\ 

\end{tabular}. \]

Let us compute the ideal $I$ defining the stabilizer $\St_G$ of the groupoid $G \times_X Y \toto Y$. It is the ideal generated in $A[G]$ by the image of the ideal defining the diagonal immersion \[ Y \into Y \times_X Y. \]

The latter is generated in $A \otimes_{A_0} A$ by the elements of the form $(1 \otimes a - a \otimes 1)$, for $a \in A$. Note that we have \[ j^\sharp(1 \otimes a - a \otimes 1) = \displaystyle\sum_{g \in G} (g(a)-a)e_g  \] 

since in $A[G]$ we have $1 = \displaystyle\sum_{g \in G} e_g$. These expressions generate the ideal $I$.

The augmentation ideal of $\St_G$ is generated by the images in $A[G]/I$ of the $e_g$ with $g \neq 1$.

Observe that, since $e_g e_{g'}=\delta_{g,g'}$, if $t = \displaystyle\sum_{g \in G} t_ge_g \in A[G]$ and $u=\displaystyle\sum_{g \in G } (g(a)-a)e_g$ for some $a \in A $ then \[tu = \displaystyle\sum_{g \in G} t_g(g(a)-a)e_g. \]

We thus have an isomorphism of $A$-algebras \[ \cO_{\St_G} \simeq \displaystyle\bigoplus_{g \in G} A/I_g, \] where $I_g$ is the ideal generated in $A$ by the expressions $(g(a)-a)$, $a \in A$. 

It follows that we have the isomorphism of $A$-modules \[ \fm_G \simeq \displaystyle\bigoplus_{g \neq 1} A/I_g. \] 

By assumption, the residue field extension $k(A)/k(A_0)$ is separable. Hence by \cite[III, \S 6, prop.12]{Se}, $A$ is a monogenic $A_0$-algebra. Let $x$ be a generator and $v$ be the valuation in $A$. For all $g \in G$ we have \[v(I_g) = v(g(x)-x) := i_G(g). \] With these notations we thus have \[ \fm_G \simeq \displaystyle\sum_{g \in G} A/ \pi^{i_G(g)}, \] where $\pi$ is a uniformizer of $A$. 

We obtain a free resolution of the $A$-module $\fm_G$ of the following form : \[ 0 \too A^{\oplus \vert G \vert -1} \stackrel{M}{\too} A^{\oplus \vert G \vert -1} \too \fm_G \too 0, \] where $M$ is a diagonal matrix of size $\vert G \vert -1$ whose diagonal entries are the elements $\pi^{i_G(g)}$ for $g \neq 1$. Its determinant has valuation $ \displaystyle\sum_{g \neq 1} i_G(g)$. 

On the other hand, we know from \cite[IV, \S 1, prop.4]{Se} that \[\displaystyle\sum_{g \neq 1} i_G(g) = v( \fD_{A/A_0} ), \] where $\fD_{A/A_0}$ stands for the different of the ring extension $A/A_0$. The latter is the annihilator of the module $\Omega^1_{A/A_0}$.

Thus we see that the multiplicities of the divisors $R_\cG$ and $R_{Y/X}$ are equal at every codimension $1$ points. Hence they are equal.

\end{demo}

\begin{rema}

The term "generalized covering" that we use to refer to our definition \ref{def} is abusive because it is not clear to us how to include the generically \'etale morphisms that do not arise from group actions. More precisely, if $Y \too X$ is a finite locally free generically \'etale morphism we do not know what groupoid to attach to it in order to make it a generalized covering in the sense of \ref{def}. One can always consider the trivial groupoid $Y \times_X Y \toto Y$ given by the two projections but this would not be a wise choice since its stabilizer is always trivial, even if $Y \too X$ is not \'etale everywhere.  

On the other hand, if $Y \too X$ as above is given by the quotient of $Y$ by the action of a finite group, then it becomes a generalized covering  in the sense of \ref{def} when endowed with the action groupoid.

\end{rema}

\subsection{Devissage of the ramification divisor.}

We now tackle the problem of performing d\'evissage of generalized coverings. Let $X$ be an $S$-scheme and $(f : Z \too X , \cG \toto Z)$ be a generalized covering of $X$. Suppose given a subgroupoid $\cH \toto Z$ of $\cG \toto Z$. Since $\cG$ acts generically freely on $Z$, so does $\cH$. Let $Y = Z/\cH$ be the quotient scheme of $Z$ by $\cH$. We then have a covering $(g : Z \too Y ,\cH \toto Z)$ of $Y$. Since obviously the morphism $f : Z \too X$ is $\cH$-invariant, we have a factorisation  

\[ \xymatrix{ Z \ar[rr]^g \ar[rd]_f && Y \ar[ld]^h \\ & X } .\]

We would like to add some structure on the morphism $h$ in order to make it into a covering of $X$.
The scheme $Y$ should be endowed with the action of the quotient groupoid of $\cG$ by $\cH$. By quotient groupoid we mean a groupoid $\cQ \toto  Y$ acting on $Y$ such that every groupoid morphism $(\cG \toto Z) \too (\cT \toto T)$ that contains $\cH$ in its kernel factors through $\cQ \toto Y$. In case the groupoid $\cG \toto Z$ is given by the action of a finite group scheme $G$ on $Z$ and $H \lhd G$ is a normal subgroup, it is easy to check that the quotient groupoid is just the residual action groupoid $G/H \times_S Y \toto Y$ of $G/H$ on $Y$.

In general, since we want to have a groupoid with source and target defined in $Y$, it is natural to define $\cQ$ as the quotient scheme of $\cG$ by the action of $\cH^2$ by pre- and post-composition, ie we define $\cQ$ as the quotient of the groupoid \[ (\cH \times_Z \cH) \times_{(s,s),Z \times_S Z, (t,s)} \cG \toto \cG \] whose arrows are of the form \[ (\varphi,\psi,g) : g \too \varphi g  \psi^{-1}. \]
For short, we denote it by $\cG' \toto \cG$. The compositions $\cG \toto Z \too Y$ are invariant under the action of the above groupoid so we get maps $\sigma, \tau  : \cQ \toto Y$ which will be the source and target of the groupoid we wish to define. It is however not obvious to us how to define the composition of arrows in $\cQ$, ie how to fill the diagram \[ \xymatrix{\cG \times_{s,Z,t} \cG \ar[r]^-c \ar[d] & \cG \ar[d] \\ \cQ \times_{\sigma,Y,\tau} \cQ & \cQ } \]  
in a systematic way. This is because, since the actions involved to construct the quotients are a priori not free, we do not know their points and hence cannot just lift points in $\cQ \times_{\sigma,Y,\tau} \cQ$, compose their lifts in $\cG \times_{s,Z,t} \cG$ and send the result back in $\cQ$. This separate problem will be the subject of a subsequent paper. 

\medskip

 In the sequel of this section we assume that such a quotient groupoid has been constructed and we investigate the behaviour of the ramification divisor under such a d\'evissage. We fix an $S$-scheme $X$, a covering $(f : Z \too X, \cG \toto Z)$ of $X$ and a subgroupoid $\cH \into \cG$. We denote the quotient of $Z$ by $\cH$ by $g : Z \too Y$. We let $\cQ$ be the quotient of the groupoid $\cG' \toto \cG$ defined above.

First we must check that we have the following lemma :

\begin{lemm}

With the above notations, suppose that a quotient groupoid $\cQ \toto Y$ has been constructed. If $ g : Z \too Y$ is flat then

\begin{itemize}

\item[(i)] $(g : Z \too Y, \cH \toto Y)$ is a covering of $Z$.

\item[(ii)] $(h : Y \too X, \cQ \toto Y)$ is a covering of $X$.

\end{itemize}

\end{lemm}

\begin{demo}

Only $(ii)$ needs a proof. If $\cG$ acts freely on $Z$, so does $\cH$ and then $Y$ represents the fppf quotient sheaf $T \mapsto Z(T)/\cH(T)$. In the same way $\cQ$ represents the quotient sheaf $T \mapsto \cG(T)/\cH^2(T)$. We then easily see that $\cQ$ acts freely on $Y$. 
Now if $U \subset Z$ is a dense open subscheme on which $\cG$ acts freely, by faithfull flatness of $g : Z \too Y$, its image $V \subset Y$ is a dense open subscheme of $Y$, on which $\cQ$ acts freely by the above discussion. Furthermore the fiberwise criterion for flatness shows that $h$ is flat, since both $f$ and $g$ are. 
\end{demo}

The quotient morphisms $p : \cG \too \cQ$ and $g : Z \too Y$ induce a groupoid morphism which we still denote $p : (\cG \toto Z) \too (\cQ \toto Y)$. The following lemma shows that, in our situation and under flatness assumptions its kernel will be $\cH$, as expected.

\begin{prop}
\label{ker}

With the above notations, if $\cG \too \cQ$ and $Z \too Y$ are flat, we have $\ker p \simeq \cH$.
\end{prop}

\begin{demo}

Denote by $\cI_\cH \subset \cO_\cG$ the ideal sheaf defining $\cH$ in $\cG$. We want to show that it agrees with the ideal sheaf defining $\ker(p)$. This is a local question so we may assume that all schemes involved are affine and work with global sections. Let $x \in \cI_\cH \cap \cO_\cQ$. It defines a morphism $\cG \too \AA^1$, $\cH$-invariant and vanishing on $\cH$. We have the following commutative diagram : 

 \[ \xymatrix{ & Z \ar[r]^g \ar[d]_{e} & Y \ar[d]^{\bar{e}}   \\ \cH \ar[r]^i \ar[ur]^{s} & \cG \ar[d]_x \ar[r]^p  & \cQ \ar[ld]^{\bar{x}} \\ & \AA^1_S & &} \]
  
where $e$ (resp. $\bar{e}$) is the unit section of the groupoid $\cG \toto Z$ (resp. $\cQ \toto Y$) and $\bar{x}$ is the function $\cQ \too \AA^1_S$ induced by the $\cH$-invariant function $x$. 

We have $\bar{x}\circ \bar{e} \circ g = x \circ e = 0$ and since $g$ is an epimorphism we have $\bar{x} \circ \bar{e} = 0$ and thus $x \in \fm_\cQ$. Hence $\cI_\cH \cap \cO_\cQ \subset \fm_\cQ$. 

Conversely, since $x \in \fm_\cQ$ we have $\bar{x} \circ \bar{e} = 0$. Since by the above diagram we have $(\bar{x} \circ \bar{e}) \circ (g \circ s) = x \circ i$, we see that $x$ vanishes on $\cH$ and $x \in \cI_\cH \cap \cO_\cQ$. We thus have $\cI_\cH \cap \cO_\cQ = \fm_\cQ$. 

The kernel of $p$ is defined by the following fibre product : 
\[ \xymatrix{ \ker p \ar[d] \ar[r] & Y \ar[d]_{\bar{e}} \\ \cG \ar[r] & \cQ}. \]

Its structure sheaf is thus $\cO_{\cG} \otimes_{\cO_{\cQ}} \cO_Y = \cO_{\cG} / \fm_{\cQ} \cO_{\cG}$. Since $\cI_\cH \cap \cO_\cQ = \fm_\cQ$ we have $\fm_\cQ \cO_\cG \subset \cI_\cH$. We thus have a surjection $\cO_\cG / \fm_\cQ \cO_\cG \too \cO_\cG / \cI_\cH$ and a closed immersion $\cH \into \ker(p)$.

Let $U \subset Z$ be a dense open subscheme on which $\cH$ acts freely. We may assume that $U$ is saturated, ie $U = g^{-1}(g(U))$. By flatness the preimage $U'$ of $U$ in $\cG$ is also dense. Then the action of $\cG'_{\vert U'}$ on $U'$ is also free. Indeed, for $(\varphi,\psi,g) \in \cG'$ the equality $g = \varphi g \psi^{-1}$ implies that $g \psi = \varphi$ and thus that $s(\psi) = t(\psi) = s(g)$ and $s(\varphi) = t(\varphi) = t(g)$. Hence $\psi$ and $\varphi$ must be in the stabilizer of $\cH_{U}$, which is trivial. 
Let us denote $V = g(U)$. By faithful flatness of $g$, it is dense in $Y$. 
Since the formation of quotients commute with flat base change, $\cQ_{\vert V}$ (resp. $V$) represents the fppf of $U'$ (resp. $U$) by $\cG_{\vert U'}'$ (resp. $\cG_{\vert U}$).  

We can then verify on points that $\cH_{\vert U}$  is the kernel of the projection $\cG_{\vert U } \too \cQ_{\vert V}$. Indeed, if $T$ is an $S$-scheme and $t \in \cG_{\vert U} (T)$ is such that $p(t) = 1_{g(s(t))}$ there exists an fppf covering $T' \too T$ and $\varphi,\psi \in \cH(T')$ such that, restricting to $T'$, we have $1_{s(t)} = \varphi t \psi^{-1}$ and hence $t = \varphi^{-1} \psi \in \cH(T')$.

Thus we see that the immersion $\cH \too \ker(p)$ induces an isomorphism $\cH_{\vert U} \simeq \ker(p)_{\vert U}$. By assumption $\cG \too \cQ$ is flat so $\ker(p) \too Y$ is flat. Hence the preimage of $U$ in $\ker(p)$ is dense. The closed immersion $\cH \into \ker(p)$ is thus dominant, so it is an isomorphism.

\end{demo}

We will use the above proposition to relate the different stabilizers involved. We first need the following two general observations about morphisms of groupoids. 

\begin{lemm}

\label{lemm}
Let $ a : (\cA \toto A) \too (\cB \toto B)$ be a morphism of $S$-groupoids.
\begin{itemize}
\item[(i)] $a$ induces a morphism of $S$-group schemes $a' : \St_\cA \too \St_\cB$ and we have $\ker(a') = \ker(a) \times_{\cA} \St_\cA$.

\item[(ii)] let $a_0 : A \too B$ be the morphism of object schemes induced by $a$. The morphism $a$ induces a morphism of groupoids $\tilde{a} : (\cA \toto A) \too (a_0^*\cB \toto A)$ and we have $\ker(\tilde{a}) = \ker(a) \times_{\cA} \St_\cA$.

\end{itemize}

\end{lemm}

\begin{demo}

\begin{itemize}

\item[(i)] Since $a$ is groupoid morphism it maps $\St_\cA$ to $\St_\cB$ so it induces a morphism $a'$ as stated. We have $\ker(a) = B \times_\cB \cA$ so $\ker(a) \times_\cA \St_\cA = B \times_\cB \St_\cA$. Since $B \too \cB$ factors through the immersions $B \into \St_\cB \into \cB$ we see that $B \times_\cB \St_\cA = B \times_{\St_\cB} \St_\cA = \ker(a')$.

\item[(ii)] Recall that $a_0^*\cB = \cB \times_{B \times_S B} A \times_S A$  Since $a$ is a groupoid morphism the diagram \[\xymatrix{\cA \ar[r] \ar[d] & \cB \ar[d] \\ A \times_S A \ar[r] & B \times_S B } \] is commutative. We thus get a morphism   
 \[
\begin{tabular}{cccc}
 
$\tilde{a} : $ & $\cA$ & $\too$ & $f^* \cB$ \\ 
 
& $r$ & $\mapsto$ & $(f(r),s(r),t(r))$ \\ 
 
\end{tabular} \]

where $s$ and $t$ are the source and target in $\cA \toto A$. The unit section of $a_0^* \cB$ is given by 
\[ \begin{tabular}{ccc}
$A$ & $\too$ & $f^* \cB$ \\ 

$a$ & $\mapsto$ & $(\id_{f_0(a)}, a,a)$ \\ 

\end{tabular} \]
so $\ker(\tilde{a})= \ker (a) \times_\cA \St_\cA$.
\end{itemize}
\end{demo}

We can now state our result, which we view as a substitute for the first fundamental exact sequence of the sheaf of differential $1$-forms. 
\begin{theo}
\label{exactgroup}
Let $X$ be an $S$-scheme, $(Z \too X, \cG \toto Z)$ be a covering of $X$ and $\cH \into \cG$ be a subgroupoid of $\cG$. Denote by $g : Z \too Y$ the quotient of $Z$ by $\cH$. Suppose constructed the quotient groupoid $\cQ \toto Y$.

We have the following exact sequence of $Z$-group schemes : \[ 1 \too \St_\cH \too \St_\cG \stackrel{\alpha}{\too} g^*\St_\cQ, \] 

in the sense that 

\begin{itemize}

\item[(i)] $\St_\cH \too \St_\cG$ is a closed immersion.

\item[(ii)] $\ker(\alpha) \simeq \St_\cH$.
\end{itemize}

\end{theo}

\begin{demo}
The morphism $p : \cG \too \cQ$ induces a morphism $\tilde{p} : \cG \too g^*\cQ$ whose kernel is $\ker(p) \times_\cG \St_\cG$ by the lemma \ref{lemm}. Let $\alpha : \St_\cG \too g^*\St_\cQ$ be the induced morphism on stabilizers.

\begin{itemize}

\item[(i)] $\St_\cH \too \St_\cG$ is the base change of the closed immersion $\cH \into \cG$ by the diagonal $Z \into Z \times_Y Z$. Thus it is a closed immersion.

\item[(ii)] By the lemma \ref{lemm} we have $\ker(\alpha) = \ker(\tilde{p})$. By proposition \ref{ker} we have $\ker(p) \simeq \cH$. Hence $\ker(\alpha) = \ker(p) \times_\cG \St_\cG \simeq \cH \times_\cG \St_\cG = \St_\cH$.

\end{itemize}

\end{demo}

Accordingly, we obtain an exact sequence relating the augmentation ideals of the stabilizers involved.

\begin{coro}

\label{exactmodule}

With the notations and hypothesis of the previous theorem, we have the following exact sequence of $\cO_Z$-modules : 

\[ 0 \too g^*\fm_\cQ \cO_{\St_\cG} \too \fm_\cG \too \fm_\cH \too 0, \]

where $\fm_\cQ$, $\fm_\cG$ and $\fm_\cH$ respectively stand for the augmentation ideals of the group schemes $\St_\cQ$, $\St_\cG$ and $\St_\cH$.

\end{coro}

\begin{demo}

If $\cI_\cH$ is the ideal sheaf defining $\St_\cH$ in $\St_\cG$ we have the following exact sequence \[ 0 \too \cI_\cH \too \cO_{\St_\cG} \too \cO_{\St_\cH} \too 0. \] Since $\ker(\alpha) = \St_\cG \times_{g^*{\St_\cQ}} Z$, the ideal sheaf  defining $\ker(\alpha)$ is the one generated by $g^*\fm_\cQ$ in $\cO_{\St_\cG}$, namely $g^*\fm_\cQ \cO_{\St_\cG}$. By \ref{exactgroup} we have $\St_\cH \simeq \ker(\alpha)$, hence $\cI_\cH = g^*\fm_\cQ \cO_{\St_\cG}$. Thus we have the exact sequence \[ 0 \too g^*\fm_\cQ \cO_{\St_\cG} \too \cO_{\St_\cG} \too \cO_{\St_\cH} \too 0. \label{suite} \tag{1} \] 

But since the unit section splits the structure maps $\St_\cG \too Z$ and $\St_\cH \too Z$, as $\cO_Z$-module we have $\cO_{\St_\cG} = \fm_\cG \oplus \cO_Z$ and $\cO_{\St_\cH} = \fm_\cH \oplus \cO_Z$. Modding out by $\cO_Z$ in (\ref{suite}), we get the exact sequence announced.

\end{demo}

Since the length of modules is additive with respect to exact sequences, taking the associated divisors of the modules involved in \ref{exactmodule}, we get the following equality between the associated divisors on $Z$ : 

\[ \R_\cG = \R_\cH + \textup{div}(g^*\fm_\cQ \cO_{\St_\cG}). \]

\label{transtriv} In particular if the quotient groupoid $\cQ$ acts freely on $Y$, we get the equality $\R_\cG = \R_\cH$, which is a weak form of the transitivity property (\ref{trans}) for the ramification divisor in this context. 

It should be noted that the map $\alpha$ is neither flat nor dominant. Hence it is not obvious to relate the sheaf $ g^*\fm_\cQ \cO_{\St_\cG}$ to $g^*\fm_\cQ$. In general, the formula $\R_\cG = \R_\cH + g^* \R_{\cQ}$, that one might expect in analogy with the classical situation, is not true, as illustrated by the following example.

\begin{exam}

Let $n$ be a positive integer and $k$ a field of positive characteristic $p>0$. Consider the action of the group scheme $\GL_{n,k}$ of invertible $n \times n$ invertible matrices on the $n \times n$ matrices $\M_{n,k} = \Spec(k[z_{ij}])$ by left multiplication : 
\[ \begin{tabular}{ccc}
$\Gl_n \times_k \M_n$ & $\too$ & $\M_n$ \\ 

$(P,M)$ & $\mapsto$ & $PM$ \\ 

\end{tabular} \]

For all positive integer $\gamma$, let $\G_\gamma$ be the kernel of the $\gamma$-th iterated Frobenius morphism \[ F_\gamma : \GL_n \too \GL_n^{(\gamma)}. \]
We have \[ \G_\gamma = \Spec \left( \frac{k[a_{ij}, \ 1 \leq i,j \leq n]}{a_{ij}^{p^\gamma} \ i \neq j \ , \ a_{ii}^{p^\gamma} -1} \right). \] 

It is a finite flat $k$-group scheme of order $p^{\gamma n^2}$. The above action of $\GL_n$ induces an action of all the Frobenius kernels, by the same formula.

Let $0< \beta < \gamma$ be two integers. By \cite[I, \S 9.4-9.5]{Jantzen}, $\G_\beta$ is a normal subgroup of $\G_\gamma$ and $\G_\gamma/G_\beta \simeq \G_{\gamma-\beta}$. Set $Z = \M_{n,k}$, $Y = Z/\G_\beta$ and $X=Z/\G_\gamma$. For $\tau \in \{\beta,\gamma,\gamma-\beta\}$ we denote by $\cG_\tau$ the action groupoid associated with the action of $\G_\tau$ and by $\St_\tau$ its stabilizer.  

We thus have defined a covering $(Z,\cG_\gamma)$ of $X$ which we expressed as the covering $(Z,\cG_\beta)$ of $Y$ followed by the covering $(Y,\cG_{\gamma-\beta})$ of $X$.

The quotient morphism $g : Z \too Y$ is easily seen to be defined by the ring map \[ \begin{tabular}{ccc}

$k[y_{ij}]$ & $\too$ & $k[z_{ij}]$ \\ 

$y_{ij}$ & $\mapsto$ & $z_{ij}^{p^\beta}$ \\ 

\end{tabular}. \]

For $\tau \in \{\beta,\gamma \}$ one can show that the stabilizer of the corresponding action is given by \[ \cO_{\St_\tau} = \frac{k[z_{ij}][a_1,\dots,a_n]}{a_1^{p^\tau},\dots,a_n^{p^\tau}, \Delta a_1,\dots,\Delta a_n} \] 

Similarly, we have $\cO_{\St_{\gamma-\beta}} = \frac{k[y_{ij}][b_1,\dots,b_n]}{b_1^{\gamma-\beta},\dots,b_n^{\gamma-\beta}, \Delta'b_1,\dots,\Delta'b_n}$, where $\Delta$ (resp. $\Delta'$) stands for the determinant polynomial in the variables $z_{ij}$ (resp. $y_{ij}$). 

The corresponding ramification divisors are thus : $\R_\gamma = (p^{n\gamma}-1)[\Delta]$, $\R_\beta = (p^{n\beta}-1)[\Delta]$ and $\R_{\gamma-\beta}= (p^{n(\gamma-\beta)}-1)[\Delta']$.
Since the determinant $\Delta'$ in the variables $y_{ij}$ is mapped to $\Delta^{p^\beta}$ we have $g^*\R_{\gamma - \beta} = p^{\beta}(p^{n(\gamma - \beta)} -1)[\Delta]$.
Thus we see that $\R_\gamma \neq \R_\beta + g^*\R_{\gamma-\beta}$.

\end{exam}

\section{Generalized coverings given by diagonalizable group actions.}

\subsection{Diagonalizable group schemes and their actions.}

We briefly recall some definitions and facts concerning diagonalizable group schemes and their actions which will be useful for us. We refer to \cite[Exp. VIII]{SGA3} for details and proofs.

\begin{defi}

An $S$-group scheme $G$ is said to be diagonalizable if it is isomorphic to the character group scheme of a constant group, ie if there exists an abstract abelian group $M$ and an isomorphism of $S$-group schemes $G \simeq \Hom_{Grp/S}(M_S, \GG_{m,S})$, where $M_S$ is the constant $S$-group scheme defined by $M$. This is equivalent to the existence of an isomorphism of $\cO_S$-Hopf algebras $\cO_G \simeq \cO_S[M]$.

If $M$ is an abstract group we often denote by $\D(M)$ the $S$-group scheme associated to it. We obtain a contravariant functor $M \mapsto \D(M)$ from abstract groups to diagonalizable $S$-group schemes.

\end{defi}

We have the following lemma, from \cite[Exp VIII, prop.2.1]{SGA3} :

\begin{lemm}

The group scheme $\D(M) \too S$ is smooth if and only $M$ is of finite type and the order of its torsion subgroup is prime to all the residue characteristics of $S$. 
A diagonalizable group scheme is always faithfully flat and affine over $S$.

\end{lemm}

For an integer $n \geq 2$ we denote by $\mu_{n,S}$ the group scheme $\D(\ZZ/n\ZZ)$. According to the lemma above, it is \'etale if and only if $n$ is prime to the characteristics of all residue fields of $S$.

It turns out that actions of diagonalizable group schemes are easy to describe in terms of graded algebras. Let us recall the following definition :

\begin{defi}

Let $M$ be an abstract abelian group. An $\cO_S$-algebra $\cA$ is said to be $M$-graded if it has a decomposition as $\cO_S$-module

\[ \cA = \displaystyle\bigoplus_{m \in M} \cA_m, \] 
where \begin{itemize}
\item[•] $\cA_0$ is a sub-$\cO_S$-algebra of $\cA$
\item[•] For all $(m,n) \in M^2$, $\cA_m \cA_n \subset \cA_{m+n}$

\end{itemize}

\end{defi}

We have the following proposition, from \cite[Exp I, 4.7.3]{SGA3} and \cite[Exp VIII, Prop 4.1-4.6]{sg32} : 

\begin{prop}
\label{diago}

Let $M$ be an abstract abelian group. The functor $\cA \mapsto \Spec(\cA)$ induces an anti-equivalence between the category of $M$-graded quasi-coherent $\cO_S$-algebras and the category of affine $S$-schemes with an action of $\D(M)$.

\medskip

Let $Y=\Spec(\cA)$ be an $S$-scheme with an action of $\D(M)$. Then $Y \too S$ is a $\D(M)$-torsor if and only if the following two conditions are satisfied : 

\begin{itemize}

\item[(a)] For all $m \in M$, $\cA_m$ is an invertible $\cO_S$-module.
\item[(b)] For all $(m,n) \in M^2$, the morphism $\cA_m \otimes_{\cO_S} \cA_n \too \cA_{m+n}$ induced by multiplication is an isomorphism.
\end{itemize}
Those two conditions are in turn equivalent to the following : 

\begin{itemize}

\item[(a')] The morphism $\cO_S \too \cA_0$ is an isomorphism.
\item[(b')] For all $m \in M$, $\cA_m \cA_{-m} = \cA_0$.

\end{itemize}

\end{prop}

We analyse the structure of covegins given by actions of diagonalizable group schemes. We are of course interested in the non-\'etale case.

\subsection{Local structure of $\D(M)$-coverings.}

Let us fix an $S$-scheme $X$. Let $G = \D(M)$ be a finite diagonalizable group scheme acting on an $S$-scheme $Y$. Suppose given a morphism $f : Y \too X$ such that $(Y,G)$ is a covering of $X$. Then $Y$ is affine over $X$. As recalled above, there exists an $M$-graded $\cO_X$-algebra $\cA$ such that $Y = \Spec(\cA)$. 

The action of $G$ on $Y$ is given by the map 

\[\begin{tabular}{ccc}

$\cA$ & $\too$ & $\cA \otimes_{\cO_S} \cO_S[M]$ \\ 

$a = \displaystyle\sum_{m \in M}a_m$ & $\mapsto$ & $\displaystyle\sum_{m \in M} a_m \otimes X^m$ \\ 

\end{tabular}, \]

where we used the notation $X^m$ to denote the generator of $\cO_S[M]$ coreesponding to $m \in M$.

The sub-algebra of invariants is $\cA_0$, so $X = Y/G = \Spec(\cA_0)$.
  
By definition of a covering, there exists a schematically dense open subscheme $V \subset Y$ on which $G$ acts freely. Replacing $V$ by its $G$-orbit if necessary, we may assume that $V$ is $G$-stable. The morphism 
\[ \begin{tabular}{cccc}

$j_X$ :  & $G \times_S Y$ & $\too$ & $Y \times_X Y$ \\ 

 & $(g,y)$ & $\mapsto$ & $(y,g.y)$ \\ 

\end{tabular}   \]

then induces an isomorphism $G \times_S V \simeq V \times_X V$.

Since $G \too S$ is flat, $G \times_S V$ is schematically dense in $G \times_S Y$, so $j_X$ is schematically dominant. The map \[ j_X^\sharp : \cO_Y \otimes_{\cO_X} \cO_Y \too (j_X)_* \cO_G \otimes_{\cO_S} \cO_Y \]

is thus injective.
Since $f : Y \too X$ is locally free, each $\cA_m$ is locally free. One can thus cover $X = \cup_i U_i$ by open affines such that the restriction of each of the $\cA_i$ are free in restriction to the open affines $f^{-1}(U_i) \subset Y$. 
Hence to investigate the structure of the covering $(Y,G)$, we may assume that $S$, $X$ and $Y$ are affine schemes and that $\cO_Y := \cA$ is free (necessarily of finite rank) over $\cO_X$. Let us denote $Y = \Spec(A)$, $X = \Spec(A_0)$ and $S = \Spec(B)$. Up to further localization in $X$ if necessary, we may also assume that each of the pieces $A_m$ of the $M$-grading of $A$ is free. To compute the rank of $A_m$ we may restrict to the image in $X$ of a dense open subscheme of $Y$ on which the action of $G$ is free. Then condition $(b')$ of \ref{diago} shows that the rank of $A_m$ as an $A_0$-module is $1$. \label{rank}

We then have the following result, giving the local structure of covering under diagonalizable groups :

\begin{theo}

\label{struct}

With the previous conventions and notations, there exists a basis $(e_m)_{m \in M}$ of $A$ as an $A_0$-module with $e_0 =1$ and non zero-divisors  $(\alpha_{m,n})_{m,n \in M}$ in $A_0$ with $\alpha_{0,n} = \alpha_{m,0} = 1$, $\alpha_{m,n}=\alpha_{n,m}$ and \[ \forall l,m,n \in M, \ \alpha_{l,m}\alpha_{l+m,n} = \alpha_{m,n}\alpha_{l,m+n}, \]

such that the following holds :

\begin{itemize}

\item[(i)] For all $m \in M$, $A_m = A_0 e_m$ 
\item[(ii)] For all $(m,n) \in M^2$, $e_me_n = \alpha_{m,n}e_{m+n}$

\end{itemize}

Furthermore, if $M=\ZZ/p^n \ZZ$, the $\alpha_{i,j}$ are determined by the $\alpha_{i,1}$ More precisely, let $i \in \ZZ / p^n \ZZ$. Denote by $s(i)$ the unique integer in $\{0,\dots,p^n-1\}$ whose class modulo $p^n$ is $i$. Define the map 

\[
 \begin{tabular}{cccc}

$\sigma$ :  & $\ZZ / p^n \ZZ \times \ZZ / p^n \ZZ$ & $\too$ & $\ZZ$ \\ 

 & $(i,j)$ & $\mapsto$ & $\frac{1}{p^n}(s(i) + s(j) - s(i+j))$ \\ 

\end{tabular} \]

and set $\beta_0=1$, $\beta_{i+1} = \alpha_{0,1} \dots \alpha_{i,1}$ for all $i \in \ZZ/ p^n \ZZ$ such that $s(i) \neq p^n-1$ and $f = \alpha_{0,1} \dots \alpha_{p^n-1,1}$. We then have, for all $i,j \in \ZZ / p^n \ZZ$, 

\begin{equation} \label{equat} \alpha_{i,j} = \beta_{i+1}\beta_i^{-1}\beta_j^{-1}f^{\sigma_{i,j}}. \end{equation}

Conversely, for all $(p^n-1)$-tuple $(\alpha_{i,1})_{i \in \{1,\dots,p^n-1 \} }$ of non zero-divisors in $A_0$ we get a $\mu_{{p^n},S}$-covering in the following way : 

\begin{itemize}

\item Set $A = A_0^{\oplus p^n}$, label each copy of $A_0$ by an index $i \in \ZZ / p^n \ZZ$ and set $e_i = (\delta_{ij})_{j \in \ZZ / p^n \ZZ}$.

\item For all $i,j \in \ZZ/ p^n \ZZ$, define $\alpha_{0,1}=1$ and $\alpha_{i,j}$ according to the formula (\ref{equat}).

\item Give $A$ the structure of a $\ZZ / p^n \ZZ$-graded $A_0$-algebra by setting, for all $i,j \in \ZZ / p^n \ZZ$, \[e_ie_j = \alpha_{i,j} e_{i+j}. \]

\end{itemize}

Then $\Spec(A)$ is a $\mu_{p^n,S}$-covering of $\Spec(A_0)$.

\end{theo}

\begin{demo}

By the preceding remarks, each of the $A_m$ is free of rank $1$ over $A_0$, for which we denote $e_m$ a generator, with the convention that $e_0=1$. Since for all $(m,n) \in M^2$ we have $A_mA_n \subset A_{m+n}$, if we denote by $\alpha_{m,n}$ the determinant of the multiplication \[ A_m \otimes_{A_0} A_n \too A_{m+n} \] we have $e_me_n = \alpha_{m,n} e_{m+n}$.

Let us note that commutativity and associativity of the multiplication in $A$ imply the following relations in $A_0$ : 

\begin{itemize}

\item[-] For all $(m,n) \in M^2$, $\alpha_{m,n} = \alpha_{n,m}$.
\item[-] For all $(l,m,n) \in M^3$, $\alpha_{l,m} \alpha_{l+m,n} = \alpha_{m,n}\alpha_{l,m+n}$.

\end{itemize}

Since the morphism $j_X^\sharp : A \otimes_{A_0} A \too A \otimes_B B[M]$ is injective between these two free $A_0$-modules of rank $\vert M \vert^2$, its determinant is a non zero-divisor. Let us compute it on the basis $(e_m \otimes e_n)$ in the source and $(e_k \otimes X^l)$ in the target. We index the matrix by $M^2$.  
We have $j_X^\sharp(e_m \otimes e_n) = e_m(e_n \otimes X^n) = \alpha_{m,n}e_{m+n} \otimes X^n$. The matrix of $j_X^\sharp$ in these basis is thus monomial~: its coefficient of index $((k,l),(m,n))$ is zero if $(k,l) \neq (m+n,n)$ and $\alpha_{m,n}$ otherwise. We thus have \[ \det(j_X^\sharp)= \varepsilon(\tau) \displaystyle\prod_{(m,n) \in M^2} \alpha_{m,n}, \] where $\varepsilon(\tau)$ is the signature of the associated permutation. We conclude that the $\alpha_{m,n}$ are all non zero-divisors.  \label{nzero}

Let us specify to the case where $M = \ZZ / p^n \ZZ$, ie $G = \mu_{p^n,S}$.

Inside the localization of $A_0$ in the multiplicative subset of non zero-divisors we consider the multiplicative subgroup generated by the $\alpha_{i,j}$, which we denote by $N$. We consider it as a trivial $\ZZ / p^n \ZZ$-module.

By the discussion above we have, for all $i,j,k$ in $\ZZ / p^n \ZZ$, $\alpha_{i,j} \alpha_{i+j,k}\alpha_{j,k}^{-1} \alpha_{j+k,i}^{-1} = 1$. Hence the family $(\alpha_{i,j})$ defines a $2$-cocycle of $\ZZ / p^n \ZZ$ with values in $N$.
Every element $f \in N$ determines a $2$-cocycle $(f^{\sigma_{i,j}})$ in such a way that if $f=g^{p^n}$ is a $p^n$-th power then $(f^{\sigma_{i,j}})$ is the coboundary induced by the cochain $(g^{-s(i)})$. By \cite[VIII, \S 4]{Se} we have $H^2(\ZZ/ p^n \ZZ,N) = N/N^{p^n}$. Hence there exists $f \in N$ and a coboundary $\beta : \ZZ/ p^n \ZZ \too N$ such that, for all $i$ and $j$ in $\ZZ / p^n \ZZ$, \[ \alpha_{i,j}=\beta_{i+j} \beta_i^{-1}\beta_j^{-1} f^{\sigma_{i,j}}. \]

Note that the pair $(\beta,f)$ is not unique : we still obtain the cocycle $\alpha_{i,j}$ if we replace $(\beta,f)$ by $(\{f'^{s(i)}\beta_i \}, f'^{p^n}f)$ for any $f' \in N$. In particular, multiplying by $\beta_1^{-1}$ if necessary, we may assume that $\beta_1 = 1$. Let us fix $i \in \ZZ / p^n \ZZ$. The equation (\ref{equat}) with $j=0$ shows that $\beta_0 =1$; with $j=1$ we see that $\beta_{i+1} = \alpha_{i,1}\beta_if^{-\sigma_{i,1}}$. We then distinguish between two cases : 

\begin{itemize}

\item If $i \neq p^n-1$ then $\sigma_{i,1}=0$ and $\beta_{i+1} = \alpha_{i,1} \beta_i$. Thus by induction we get \[ \beta_{i+1} = \alpha_{i,1}...\alpha_{1,1}. \]

\item If $i=p^n -1$ then $\sigma_{i,1} = 1$ and since $\beta_0 = 1$ we get $1 = \alpha_{p^n-1,1} \beta_{p^n -1}f^{-1}$ and thus \[ f = \displaystyle\prod_{l \in \ZZ /p^n \ZZ} \alpha_{l,1} \]

\end{itemize}

\end{demo}

\subsection{Ramification of $\D(M)$-coverings.}

We now compute the ramification divisor of a $D(M)$-covering for a finite abelian $p$-group $M$.
 \begin{conv}
 
Our proofs in this section rely on computations of multiplicities of Weil divisors at codimension $1$ points. We thus need these points to be regular. Hence, from now on and until the end of this article, we make the additional assumption that the schemes involved in a covering are normal.

More precisely, if $X$ is an $S$-scheme and $(Y,\cG)$ is a covering of $X$ then $Y$ will always be assumed to be normal. Note that, since normality is preserved by taking invariant rings, this implies that $X$ is itself normal.  
 
 \end{conv}

Let us fix an $S$-scheme $X$, a finite abelian $p$-group $M$ and let $G = \D(M)$ be the corresponding diagonalizable group scheme. Let $(Y,G)$ be a $G$-covering of $X$. We assume that $Y$ is normal and we compute the multiplicities of the ramification divisor defined in \ref{newdiv} at each codimension $1$ point of $Y$.

Let $y \in Y$ be such a point. By assumption, the local ring $\cO_{Y,y}$ is a discrete valuation ring which we denote $A$, with valuation $v$. We have a graduation of type $M$ \[ A = \displaystyle\bigoplus_{m \in M} A_m \] induced by the action of $G$. Since $(Y,G)$ is a covering, by the discussion \ref{rank} above, each $A_m$ is a free $A_0$-module of rank $1$ for which we denote $e_m$ a generator. We choose $e_0 = 1$ and define $Y_y = \Spec(A)$.   

By definition, the formation of the stabilizer group scheme commutes with base change. Hence the stabilizer $\St_{G,y}$ of the covering $(Y,G)$ is given at $y$ by the fibre product 

\[ \xymatrix{ \St_{G,y} \ar[r] \ar[d]_\sigma & G \times_S Y_y \ar[d]^{j_y} \\
Y_y \ar[r] & Y_y \times_X Y_y }. \]

The diagonal immersion of $Y_y$ is given by the ring morphism \[ \begin{tabular}{ccc}

$A \otimes A$ & $\too$ & $A$ \\ 

$a \otimes b$ & $\mapsto$ & $ab$ \\ 
\end{tabular} \]
  
and the morphism $j_y$ is defined by the map 
\[ \begin{tabular}{ccc}

$A \otimes A$ & $\too$ & $A[M]$ \\  
$a \otimes b$ & $\mapsto$ & $\displaystyle\sum_{m \in M} ba_mX^m$ \\ 
 
\end{tabular} .\]

We thus have  \[ \cO_{\St_{G,y}} = A[M] \otimes_{A\otimes A} A  = \frac{A[M]}{ \left( \displaystyle\sum_{m \in M} ba_mX^m - ab, \ a,b \in A \right)}. \]
 
It is easily seen that the ideal defining $\St_{G,y}$ in $G \times_S Y_y$ is also generated by the elements $e_m(X^m-1)$ for $m \in M$. Hence we have 
 \[ \cO_{\St_{G,y}}  = \frac{A[M]}{e_m(X^m - 1), \ m \in M}. \]

Its augmentation is generated by the images in $\cO_{\St_{G,y}}$ of the elements $X^m - 1$ for $m \neq 0$, ie we have  \[ \fm_{\St_{G,y}} = \displaystyle\sum_{m \in M \setminus \{0\}} (X^m -1) \cO_{\St_{G,y}}. \]

Let us remark that the elements $X^m -1$ for $m \in M$ also generate the algebra $A[M]$. We thus obtain a surjection of $A$-modules \[ \varphi : A^{\oplus \vert M \vert -1} \too \sigma_* \fm_{\St_{G,y}} \]

sending a basis to the basis of the $X^m-1$ for $m \in M$. We aim at describing its kernel. 

First note that if $e_m$ is invertible in $A$ then $X^m = 1$ in $\cO_{\St_{G,y}}$. Note furthermore that the equality $e_me_n=\alpha_{m,n}e_{m+n}$ implies that if $e_m$ and $e_n$ are invertible then so is $e_{m+n}$. Thus the set \[ N := \{ n \in M \ \vert \ e_n \in A^\times \} \] is a subgroup of $M$.
Note also that if $e_m$ is invertible, its inverse must lie in $A_{-m}$. Hence we find that $e_m$ is invertible if and only if $\alpha_{m,-m}$ is (in which case $e_m^{-1} = \alpha_{m,-m}^{-1} e_{-m}$), so we also have $N = \{ n \in M \ \vert \ \alpha_{n,-n} \in A_0^{\times} \}$.

We then have the following proposition :

\begin{prop}
\label{ram}

If $N = \{0 \}$ there exists $d \in M$ such that : 

\begin{itemize}
\item[(i)] $d$ generates $M$ as an abelian group, $e_d$ is an uniformizer of $A$ and generates $A$ as an $A_0$-algebra.

\item[(ii)] The kernel of the surjection $\varphi$ is equal to the submodule $(e_dA)^{\oplus \vert M \vert -1}$ of $A^{\oplus \vert M \vert -1}$.
\end{itemize}

\end{prop}

\begin{demo}
Suppose $N = \{0\}$.

First, note that the valuations of each of the $e_m$ are disjoint. Indeed, let $(m,n) \in M^2$ such that $v(e_m)=v(e_n)$. There exists an invertible element $a \in A^\times$ such that $e_n=ae_m$. Write $a = \displaystyle\sum_{k \in M} a_ke_k$ on the basis $(e_i)$ of $A$ as an $A_0$-module. We then have \[ ae_m =\displaystyle\sum_{k \in M} a_k \alpha_{k,m} e_{m+k} = e_n. \]

Hence $a_k \alpha_{k,m} = 0$ for $k \neq n-m$ and since none of the $\alpha_{i,j}$ is a zero divisor by \ref{nzero}, we have $a_k = 0$ for $k \neq n-m$ and thus $a = a_{n-m}e_{n-m} \in A_{n-m}$. But $a \in A^\times$ then implies $e_{n-m} \in A^\times$. Since $N = \{ 0 \}$ by assumption we must have $n=m$.

We then claim that, for all $m \in M$, we have $v(e_m) \leq \vert M \vert -1$. To prove this, suppose there exists some $n \in M$ such that $v(e_n) \geq \vert M \vert$. We can then write $e_n = \pi^{ \vert M \vert} b$ for some $b \in A$, where $\pi$ is a uniformizer of $A$. But for all $x \in A$ we have $x^{ \vert M \vert} \in A_0$. Indeed, if $x = \displaystyle\sum_{m \in M} x_m$ maps to $\displaystyle\sum_{m \in M} x_m X^m $ via the coaction, since $\vert M \vert$ is a $p$-th power, $x^{\vert M \vert}$ maps to  $\displaystyle\sum_{m \in M} x_m^{\vert M \vert} X^{\vert M \vert m}$ and since $\vert M \vert$ annihilates $M$ we have $X^{\vert M \vert m}=1$ for all $m \in M$. Hence $x^{\vert M \vert}$ is invariant, ie in $A_0$. 
If we write $b = \displaystyle\sum_{k \in M}b_ke_k$ on the basis $(e_i)$ we have \[ e_n = \displaystyle\sum_{k \in M} \pi^{\vert M \vert} b_ke_k. \] By the preceding remark we have $\pi^{\vert M \vert} b_k \in A_0$ for all $k \in M$. Thus $b_k =0$ if $k \neq n$ and $\pi^{ \vert M \vert }b_n = 1$ which is absurd since $\pi$ is a uniformizer, in particular non-invertible.

Hence the valuations of the $e_m$ are all distinct and lower than $\vert M \vert -1$. We conclude that for all $i \in \{0,..., \vert M \vert -1 \}$ there exists some $m \in M$ such that $v(e_m) = i$. In particular there exists $d \in M$ such that $v(e_d)=1$, ie such that $e_d$ is a uniformizer of $A$. If $m \in M$ we can then write $e_m = a e_d^{v(e_m)}$ for some $a \in A^\times$. On the other hand we know that there exists $\gamma_m \in A_0$ such that $e_m^{v(e_m)d} = \gamma_m e_{v(e_m)d}$, namely $\gamma_m = \alpha_{m,m}\alpha_{m,2m}\dots \alpha_{m,v(e_m)d-1}$. Let us write $a = \displaystyle\sum_{k \in M} a_k e_k$. We then find \[ e_m = \displaystyle\sum_{k \in M} \gamma_m a_k \alpha_{k,v(e_m)d} e_{k +v(e_m)d}. \]

Thus we must have $a_k = 0$ if $k \neq m - v(e_m)d$ and $a = a_{m-v(e_m)d}e_{m-v(e_m)d} \in A^{\times}$. Since $N = \{0 \}$ we must have $m=v(e_m)d$ and $a=a_{m-v(e_m)} \in A_0^\times$. Thus $d$ generates $M$ as an abelian group and $e_d$ generates $A$ as an $A_0$-algebra.

In $\cO_{\St_{G,y}}$ we thus have $e_d(X^m-1) = 0$ for all $m \in M$ and hence $(e_dA)^{\oplus \vert M \vert -1} \subset \ker(\varphi)$. Conversely, if $x=(x_k)_{k \in M^*} \in A^{\oplus \vert M \vert -1}$ is such that $\varphi(x) = 0$ let us write $x_k = y_k + z_k$, with $y_k \in e_dA$ and $z_k \in A^\times$. Then $\varphi(x) = \displaystyle\sum_{k \in M} z_k(X^k-1)$ and since the $X^k-1$ form a basis of $k_A[M]$, where $k_A$ is the residue field of $A$, we must have $z_k = 0$ for all $k$ and thus $x \in (e_dA)^{\oplus \vert M \vert -1 }$. Hence we have $\ker(\varphi) = (e_dA)^{\oplus \vert M \vert -1}$.

\end{demo}

\begin{defi}

\label{totram}
With the preceding notations, if $N=\{0 \}$ we say that the covering $(Y,G)$ is totally ramified at $y$.

\end{defi}

The last proposition allows for the computation of the ramification divisor of a totally ramified covering.

\begin{coro}

If $(Y,G)$ is totally ramified at $y$, its ramification divisor has multiplicity $\vert G \vert - 1$ at~$y$. 

\end{coro}

\begin{demo}

Let us keep the notations of \ref{ram}. This proposition shows that we have the exact sequence of $A$-modules  \[ 0 \too A^{\oplus |M|-1} \stackrel{\times e_d}{\too} A^{\oplus |M|-1} \stackrel{\varphi}{\too} \sigma_* \fm_{\St_{G,y}} \too 0, \]

where the first arrow is the multiplication of each coordinates by $e_d$. Its determinant is $e_d^{\vert M \vert - 1}$, so has valuation $\vert M \vert -1 = \vert G \vert -1$.

\end{demo}

When $N \neq \{0\}$ the above calculations still allow us to compute the multiplicities of the ramification divisor. We just have to perform a d\'evissage to reduce to the totally ramified case.

\begin{theo}

Let $X$ be a scheme and $(Y,G)$ be a covering of $X$ given by the action of diagonalizable group scheme $G=\D(M)$ on the normal scheme $Y$. We denote by $\R_G$ its ramification divisor. For every codimension $1$ point of $Y$ there exists a maximal subgroup $H_y = \D(M/N_y)$ of $G$ such that : 

\begin{itemize}

\item[(i)] The covering $(\Spec(\cO_{Y,y}), H_y)$ of \ $\Spec(\cO_{Y,y})$ induced by the action of $G$ is totally ramified. 
\item[(ii)] The residual covering $\Spec(\cO_{Y,y})/H_y \too \Spec(\cO_{Y,y})/G$ is a $G/H_y$-torsor. 
\item[(iii)] The multiplicity of $\R_G$ at $y$ is $\vert H_y \vert -1 = \vert M/N_y\vert -1$.

\label{hy}

\end{itemize}

\end{theo}
\begin{demo}

Let us denote $\cO_{Y,y} = A$ and $Y_y = \Spec(A)$. The action of $G$ on $Y_y$ is given by a graduation \[A = \displaystyle\bigoplus_{m \in M} A_m \] of type $M$ on $A$. As before set $N_y = \{ n \in M \ \vert \ A_n \otimes_{A_0} A_{-n} \simeq A_0 \}$. It is a subgroup of $M$. For each $m \in M$, let $e_m$ be a generator of $A_m$ as an $A_0$-module and let $(\alpha_{m,n})$ be the corresponding cocyle. Applying the functor $D$ to the exact sequence \[ 0 \too N_y \too M \too M/N_y \too 0 \] we get the exact sequence of group schemes \[ 1 \too H_y \too G \too G/H_y \too 1 .\]

Thus the quotient by $G$ can be factored into the quotient by $H_y$ followed by the quotient by $G/H_y$, according to the diagram 

\[ \xymatrix{ Y_y \ar[rr] \ar[rd] && Y_y/H_y \ar[ld] \\ & Y_y/G }. \]

This can be written on the type $M$ graduation of $A$ as follows : for each $\bar{m} \in M$ choose a lift $m \in M$ and write \[ A = \displaystyle\bigoplus_{\bar{m} \in M/N} \displaystyle\bigoplus_{n \in N} A_{m+n}, \]

the quotient morphisms being given by the inclusions $A_0 \into \displaystyle\bigoplus_{n \in N} A_n \into A$.

By definition for all $n \in N$ the multiplication $A_n \otimes_{A_0} A_{-n} \too A_0$ is an isomorphism so $Y_y/H_y \too Y_y/G$ is a $G/H_y$-torsor. Hence its ramification divisor is trivial and by \ref{transtriv} on a neighborhood $U$ of $y$ we have $(\R_G)_U = (\R_{H_y})_U$

We then claim that $Y_y \too Y_y/H_y$ is totally ramified. To see this, set $B_0 = \cO_{Y_y/H_y} = \displaystyle\bigoplus_{n \in N} A_0e_n$. For every residue class $\bar{m} \in M/N_y$, fix a lift $s(\bar{m}) \in M$.
Then the action of $H_y$ on $Y_y$ can be described by the graduation $A = \displaystyle\bigoplus_{\bar{m} \in M/N_y} B_0 e_{s(\bar{m})}$ and the cocycle $(\beta_{\bar{m},\bar{m'}} = \alpha_{s(\bar{m}),s(\bar{m'})})$. By definition, $\bar{m} \in M/N_y$ we have $e_{s(\bar{m})} \in A^\times$ if and only if $s(\bar{m}) \in N_y$, ie $\bar{m} = 0$. Hence $Y_y \too Y_y/H_y$ is totally ramified. By the previous corollary, the ramification divisor of $(Y,G)$ has multiplicity $\vert M/N_y \vert -1 = \vert H_y \vert -1$ at~$y$.

\end{demo}

\subsection{Relation with the fixed-point scheme.}
\label{comparaison}

Another natural object to consider for describing the ramification of a $G$-covering $(Y,G)$ is its fixed point scheme, as defined in \cite{Fogarty}. Let us briefly recall its definition. Consider the fixed point functor 

\[ \begin{tabular}{cccc}

$\textup{Fix}_G$ : & $\Sch/S$ & $\too$ & $\Ens$ \\ 

& $T$ & $\mapsto$ & $\{ t \in Y(T) \ | \ \rho \circ(\id_G \times t) = p_2 \circ (\id_g \times t) \}$ \\ 

\end{tabular}, \]

where $\rho : G \times_S Y \too Y$ is the action of $G$ on $Y$ and $p_2 : G \times_SY \too Y$ is the second projection. 
In our setting, this functor is representable by a closed subscheme  of $Y$, which we denote $Y^G$. Let $\cI_G$ be the ideal sheaf defining $Y^G$. In case $G = \D(M)$ is an infinitesimal diagonalizable group scheme, we can compute $\cI_G$. Suppose $Y = \Spec(A)$ is affine and let $A = \displaystyle\bigoplus_{m \in M} A_m$ be the graduation of type $M$ associated with the action of $G$. Let $T$ be an $S$-scheme and $t : T \too Y$ be a $T$-point of $Y$, corresponding to a ring morphism $t^\sharp : A \too \cO_T(T)$. We have $t \in \textup{Fix}_G(T)$ if and only if all $a = \displaystyle\sum_{m \in M} a_m \in A$ satisfies the equality $\displaystyle\sum_{m \in M} t^\sharp(a_m)X^m = t^\sharp(a)$. Since the variables $X^m$ are $A$-linearly independent we must have \[ t^\sharp(a_m) = \begin{cases} 0  \ \text{if} & m \neq 0  \\
a_0 \  & \text{otherwise}.
\end{cases} \]   

Such is the case if and only if $t^\sharp$ factor through the quotient of $A$ by the ideal generated by the elements $e_m$ for $m \in M \setminus \{0\}$. We thus have $\cI_G = < e_m, m \neq 0 >$.

Now suppose $y \in Y$ is a point of codimension $1$. If the covering $(Y,G)$ is totally ramified at $y$, by \ref{ram} we know that one of the $e_m$ is a uniformizer of the local ring at $y$. Thus $\cI_{G,y}$ is just the maximal ideal of $\cO_{Y,y}$. Since by \ref{ram} the multiplicity of the ramification divisor at $y$ is $\vert G \vert -1 $ we have \[ \cO_Y(-\R_G)_y = (\cI_{G,y})^{\otimes \vert G \vert -1}.\]

We thus have the following proposition : 

\begin{prop}

If $(Y,G)$ is a totally ramified covering of a scheme $X$ given by the action of a diagonalizable group scheme $G$, with ramification divisor $\R_G$, for every codimension $1$ point $y \in Y$ the ideal sheaf $\cI_G$ of the fixed point scheme verifies the relation \[ \cO_Y(-\R_G)_y = (\cI_{G,y})^{\otimes \vert G \vert -1}.\]

\end{prop}

\subsection{D\'evissage of $\D(M)$-coverings}

We now investigate the behaviour of the ramification divisors through d\'evissage in the special case of coverings given by diagonalizable group actions. Unlike the general case of \ref{transtriv}, we will see that in this situation the ramification divisor behaves like the classical one with respect to d\'evissage.

\begin{prop}
Let $X$ be an $S$-scheme and $(Z,G)$ be a covering of $X$ given by the action of an infinitesimal diagonalizable group scheme $G=\D(M)$. Suppose given a subgroup $H = \D(M/N)$ of $G$ that gives rise to a covering $(Z,H)$ of the quotient $g: Z \too Y=Z/H$.  

Denote by $\R_G$, $\R_H$ and $\R_{G/H}$ the ramification divisors respectively associated to the actions of $G$, $H$ on $Z$ and of $G/H$ on $Y$. 

As divisors on $Z$ we have \[ \R_G = \R_H + g^*\R_{G/H} .\]
\label{devissage}
\end{prop}

\begin{demo}

To show this equality of divisors, we need to show that they have the same multiplicity in each codimension $1$ point of $Z$. Since $Z$ and hence $Y$ and $X$ are normal we may assume that $Z=\Spec(A)$, $Y=\Spec(B)$ and $X=\Spec(A_0)$ are spectrums of discrete valuation rings. Let $v_A$ denote the valuation in $A$.

In view of theorem \ref{hy} we may also assume that the coverings involved are totally ramified (in the sense of definition \ref{totram}). Then by proposition \ref{ram} there exists positive integers $m \leq n$ such that $G = \mu_{p^n}$, $H=\mu_{p^m}$ and $G/H = \mu_{p^{n-m}}$. By the same proposition the stabilizers of the actions are given by $\cO_{\St_\cG}= \frac{A[t]}{t^{p^n} -1,\pi(t-1)}$, $\cO_{\St_H} = \frac{A[s]}{s^{p^m}-1,\pi(s-1)}$ and $\cO_{\St_{G/H}} = \frac{B[u]}{u^{p^{n-m}}-1, \pi'(u-1)}$, where $\pi$ (resp. $\pi'$) is a uniformizer of $A$ (resp. $B$). The graduation of type $\ZZ/p^n \ZZ$ defining the action of $G$ on $Z$ is just $A = \displaystyle\bigoplus_{0 \leq i \leq p^n-1} A_0 \pi^i$.

By definition of the residual action of $G/H$ on $Y$, we have the following commutative digram \[\xymatrix{\cO_Y \ar[r]^-{\rho_{G/H}^\sharp} \ar[d]^{g^\sharp} & \frac{\cO_Y[u]}{u^{p^{n-m}}-1} \ar[d] \\ \cO_Z \ar[r]_{\rho_G^\sharp} & \frac{\cO_Z[t]}{t^{p^n}-1}}. \]

We are going to determine the valuation in $A$ of $g^\sharp(\pi')$.
Note that we have $\rho_{G/H}^\sharp(\pi') = \pi'u$. Since the quotient map $G \too G/H$ is defined by $u \mapsto t^{p^m}$, this implies that $\rho_G^\sharp(g^\sharp(\pi')) = g^\sharp(\pi')t^{p^m}$. Hence we see that $g^\sharp(\pi') \in A_0 \pi^{p^m}$. Since $A$ is of degree $p^m$ over $B$ we cannot have $v_A(   g^\sharp(\pi')) > p^m$ so $g^\sharp(\pi')$ is of valuation $p^m$.

Now $\cO_{g^*\St_{G/H}} = \cO_{\St_{G/H}} \otimes_B A = \frac{A[u]}{u^{p^{n-m}}-1, \ g^\sharp(\pi')(u-1)}$ so the length of $g^*\fm_{G/H}$ as an $A$-module is $l_A(g^* \fm_{G/H}) = (p^{n-m}-1)p^m$. Since $l_A(\fm_G) = p^n-1$ and $l_A(\fm_H) = p^m-1$ we have $l_A(\fm_G) = l_A(\fm_H) + l_A(g^*\fm_{G/H})$.

\medskip
The divisors $\R_G$ and $\R_H + g^*\R_{G/H}$ thus have the same multiplicity in each codimension $1$ point of $Z$. Hence they are equal.

\end{demo}

We wish to relate the ramification divisor of a $\D(M)$-covering to the dualizing sheaf of its quotient morphism. We first give a criterion for the latter to be Gorenstein, in the case $\D(M) = \mu_{p^n}$.

\subsection{The Gorenstein locus of a $\mu_{p^n}$-covering}

Let us first recall the following fact about dualizing sheaves of finite morphisms. See \cite[6.4.25]{Liu} for a proof.

\begin{prop}
\label{dual}
If $f : T \too T'$ is a finite locally free morphism between locally noetherian schemes then $f$ has a dualizing sheaf given by \[ \omega_f = f^!\cO_{T'}, \] where $f^!\cO_{T'} = \cH om_{\cO_{T'}}(f_*\cO_{T'},\cO_{T'})$ is viewed as an $\cO_{T}$-module via the law $t.\theta = (x \mapsto \theta(tx))$.

\end{prop}
Recall that a morphism of schemes is said to be Gorenstein if it has a dualizing sheaf which is invertible.

Let us fix a base $S$ and an $S$-scheme $X$ endowed with a $\mu_{p^n,S}$-covering whose quotient morphism we denote by $f : Y \too X$. We wish to give a necessary and sufficient condition for the morphism $f$ to be Gorenstein in terms of the structure constants of the covering. We have the following result.

\begin{theo}
Let $f : Y \too X$ be a $\mu_{p^n,S}$-covering of $X$. Denote by \[\cA = \displaystyle\bigoplus_{i \in \ZZ / p^n \ZZ} \cA_i \] the structure sheaf of $Y$, graded by the action of $\mu_{p^n,S}$. For every $i,j \in \ZZ/p^n\ZZ$, denote by $U_{ij}$ the open subscheme of $X$ over which the multiplication \[ \cA_i \otimes_{\cO_X} \cA_j \too \cA_{i+j} \] is an isomorphism.

For each $l \in \ZZ/ p^n \ZZ$ set $U_l = \cap_{i+j=l} U_{ij}$. The open subscheme of $X$ over which the morphism $f$ is Gorenstein is the union of the $U_l$ for $l \in \ZZ / p^n \ZZ$.

In particular, $f$ is Gorenstein if and only if we have \[ X = \displaystyle\bigcup_{l \in \ZZ / p^n \ZZ} U_l.\]

\end{theo}

\begin{demo}

Let us first note that, by \label{dual}, since $f$ is finite locally free of rank $p^n$, it admits a dualizing sheaf $\omega_f = \cH om(f_*\cO_Y,\cO_X)$. As $f$ is finite and locally free, $f$ is Gorenstein if and only if, for every point $y \in Y$, $\omega_{f,y}$ is free of rank $1$. 

We may thus assume that $Y=\Spec(A)$ and $X = \Spec(A_0)$ are spectrums of local rings. We denote by $(\alpha_{ij})_{i,j \in \ZZ / p^n \ZZ}$ the cocyle with values in $A_0$ inducing the action of $\mu_{p^n}$ on $Y$. The dualizing sheaf $\omega_{Y/X}$ is then given by the $A$-module $A^*:=\Hom_{A_0}(A,A_0)$ with $A$-module law $a.\theta = (x \mapsto \theta(ax))$.  

The morphism $f$ is Gorenstein at $y$ if and only if there exists a linear form $\varphi : A \too A_0$ such that $A^*= A \varphi$. Let us note $(e_i)$ the $A_0$-basis of $A$ associated to the cocycle $(\alpha_{ij})$, ie such that $e_i e_j = \alpha_{ij} e_{i+j}$ for all $i,j \in \ZZ / p^n \ZZ$ and let $(e_i^*)$ be the dual basis. Every linear form $\theta : A \too A_0$ can be written $\theta = \sum_i \theta_i e_i^*$. Thus we see that for $\varphi \in A^*$ we have $A^* = A \varphi$ if and only if for every $j \in \ZZ/p^n\ZZ$ there exists an element $b_j = \sum_i b_{ij}e_i$ in $A$ such that $e_j^* = b_j. \varphi$. 

Observe that for all triplet $i,j,k \in \ZZ/ p^n \ZZ$ we have \[ e_i . e_j^*(e_k) = e_j^*(e_i e_k) = e_j^*(\alpha_{i,k} e_{i+k}) = \alpha_{i,k} \delta_{i+k,j}, \] where $\delta$ is the Kronecker symbol, so that \[ e_i . e_j^* = \alpha_{i,j-i}e_{j-i}^*. \]

Writing $\varphi = \sum_m \varphi_m e_m^*$ we then have \begin{equation*} \label{}
\begin{split}
b_j. \varphi & = \displaystyle\sum_{i,m \in \ZZ / p^n \ZZ}b_{ij}  \varphi_m e_i.e_m^* \\
& = \displaystyle\sum_{i,m \in \ZZ / p^n \ZZ} b_{ij} \varphi_m \alpha_{i,m-i} e_{m-i}^ *   \\
& = \displaystyle\sum_{m \in \ZZ/ p^n \ZZ} \displaystyle\sum_{k+l = m} b_{k,j} \varphi_{l+k} \alpha_{k,l} e_l^*
\end{split} 
\end{equation*}

Thus $b_j.\varphi = e_j^*$ if and only if for all $l \in \ZZ / p^n \ZZ$ we have $\displaystyle\sum_{k \in \ZZ / p^n \ZZ} b_{k,j} \varphi_{l+k}\alpha_{k,l} = \delta_{l,j}$. Hence $\varphi$ generates the $A$-module $A^*$ if and only if the matrix $M(\varphi) = (\alpha_{ij} \varphi_{i+j})_{i,j \in \ZZ / p^n \ZZ}$ is invertible. With the notations of \ref{struct} we have, in the localization of $A_0$ with respect to the multiplicative subset of non zero-divisors, 
 \[\alpha_{i,j} \varphi_{i+j} = \frac{1}{\beta_i \beta_j} \beta_{i+j} \varphi_{i+j} f^{\sigma_{i,j}}. \] Let $N(\varphi)$ be the matrix $(\beta_{i+j} \varphi_{i+j} f^{\sigma_{i,j}})_{i,j \in \ZZ / p^n \ZZ}$. We then have \[ \det(M(\varphi)) = \frac{1}{\displaystyle\prod_{i \in \ZZ / p^n \ZZ} \beta_i^2} \det(N(\varphi)) \]
 
We are going to compute $\det(N(\varphi))$. Set $\gamma_i = \beta_i \varphi_i$. We can then write \[ N(\varphi) = \begin{pmatrix}
  \gamma_0 & \gamma_1 & \gamma_2  & ... &  & \gamma_{p^n-1} \\
  \gamma_1 & \gamma_2 & \gamma_3  & ... & & \gamma_0 f  \\
  . &  &  &  &  \\
   . & . & .& .& .&. & \\
   . &. & \gamma_0f& & &  \gamma_{p^n-3} f \\
   \gamma_{p^n-1} &\gamma_0f & \gamma_1f & . .. &  & \gamma_{p^n-2}f \\
\end{pmatrix}. \]

Let us note that, if $n_{ij}$ denotes the coefficient of index $(i,j)$ in $N(\varphi)$ we have 

\[ n_{ij} = \begin{cases}
 
 \gamma_{i+j} & \mbox{if} \ s(i) + s(j) \leq p^n -1\\
\gamma_{i+j} f & \mbox{otherwise}
 
 \end{cases}
 \]

Switching the $i$-th line of the matrix $N(\varphi)$ with the $(p^n-1-i)$-th for every $i \in \ZZ/ p^n \ZZ$ we obtain the following matrix : \[N'(\varphi) = 
 \begin{pmatrix}
  \gamma_{p^n-1} & \gamma_0f & \gamma_1f  & ... &  & \gamma_{p^n-2}f \\
  \gamma_{p^n-2} & \gamma_{p^n-1} & \gamma_0f  & ... & & \gamma_{p^n-3} f  \\
  . &  &  &  &  \\
   . & . & .& .& .&. & \\
   \gamma_1 & \gamma_2 & \gamma_3& & &  \gamma_{0} f \\
   \gamma_0 &\gamma_1 & \gamma_2 & . .. &  & \gamma_{p^n-1} \\
\end{pmatrix} \]
 whose coefficients are given by \[ n'_{ij} = \begin{cases}
 
 \gamma_{j-i-1} & \mbox{if} \ s(i) \geq s(j) \\
 \gamma_{j-i-1}f & \mbox{otherwise} 
 \end{cases}
 \] and thus depend only on the difference between the line and column index.

Let $\fS_{p^n}$ be the symetric group of order $p^n$, which we view as the group of bijections of $\ZZ / p^n \ZZ$. For $k \in \ZZ / p^n \ZZ$, denote by $\tau_k$ the permutation $(i \mapsto i +k)$. Observe that for all $k,l \in \ZZ / p^n \ZZ$ we have $\tau_k \circ \tau_l = \tau_{k+l}$, so that the map  
\[ \begin{tabular}{ccc}

$\ZZ / p^n \ZZ \times \fS_{p^n}$ & $\too$ & $\fS_{p^n}$ \\ 

$(k,\sigma)$ & $\mapsto$ & $\tau_k \circ \sigma \circ \tau_k^{-1}$ \\ 

\end{tabular} \]

defines an action of $\ZZ/ p^n \ZZ$ on $\fS_{p^n}$.

If $\Omega$ stands for the sets of orbits of this action we can write $\fS_{p^n} = \displaystyle\coprod_{\omega \in \Omega} \omega$ and regroup by orbits the terms in $\det(N'(\varphi))$ so as to obtain the following expression : 
\begin{equation*} 
\begin{split}
\det(N'(\varphi))  & = \displaystyle\sum_{\sigma \in \fS_{p^n} } \varepsilon(\sigma) \displaystyle\prod_{i \in  \ZZ / p^n \ZZ} n'_{i, \sigma(i)} \\
& = \displaystyle\sum_{\omega \in \Omega} \displaystyle\sum_{ \sigma \in \omega} \varepsilon(\sigma) \displaystyle\prod_{i \in \ZZ / p^n \ZZ} n'_{i,\sigma(i)} ,
\end{split}
\end{equation*} where $\varepsilon$ stands for the signature of a permutation. Observe that it is constant on orbits. Note furthermore that if $\sigma$ and $\theta$ are in the same orbit, the sets $\{ n'_{i,\sigma(i)}, i \in \ZZ / p^n \ZZ \}$ and $\{ n'_{i, \theta(i)}, i \in \ZZ / p^n \ZZ \}$ are equal. Indeed, if there exists $k$ such that $\theta = \tau_k \circ \sigma \circ \tau_k^{-1}$ we have $n_{i,\theta(i)}'= n_{i,\sigma(i-k) +k}'$ which by definition of $N'(\varphi)$ is also equal to $n_{i-k,\sigma(i-k)}'$. Thus, if $\sigma$ and $\theta$ are in the same orbit, we have \[ \varepsilon(\theta) \displaystyle\prod_{i \in \ZZ / p^n \ZZ} n'_{i,\theta(i)} = \varepsilon(\sigma) \displaystyle\prod_{i \in \ZZ / p^n \ZZ} n'_{i, \sigma(i)}. \]
By the orbit-stabilizer theorem, the number of element in each orbit is a $p$-th power. Furthermore the orbit of a permutation $\sigma$ has only one element if and only if $\sigma \circ \tau_1 = \tau_1 \circ \sigma$ since $\tau_k = \tau_1^k$, from which we see that $\sigma(i) = \sigma(0) + i$ for all $i \in \ZZ/ p^n \ZZ$, meaning that $\sigma$ is one of the $\tau_k$. They all have signature $1$.
Hence we have \[\det(N'(\varphi)) = \displaystyle\sum_{k \in \ZZ/ p^n \ZZ} \displaystyle\prod_{i \in \ZZ / p^n \ZZ} n'_{i, \tau_k(i)}.\]
Note finally that for $k \in \ZZ/ p^n \ZZ$ we have \[ n'_{i, \tau_k(i)} = \begin{cases}

\gamma_{k-1} & \mbox{if} \ s(i) \geq s(i+k) \\
\gamma_{k-1} f & \mbox{otherwise} 
\end{cases} \] so that $\displaystyle\prod_{i \in \ZZ / p^n \ZZ} n'_{i, \tau_k(i)} = \gamma_{k-1}^{p^n} f^{p^n-1-k}$ and hence we have \[\det(N(\varphi)) = (-1)^{\frac{p^n-1}{2}} \det(N'(\varphi)) = (-1)^{\frac{p^n-1}{2}} \displaystyle\sum_{k \in \ZZ/ p^n \ZZ}\gamma_{k-1}^{p^n} f^{p^n-1-k}. \]

We now wish to come back to $M(\varphi)$. 

For all $i \in \ZZ/ p^n \ZZ$, set $c_i = \frac{ \beta_i^{p^n} f^{p^n-1-i}}{\displaystyle\prod_{j \in \ZZ / p^n \ZZ} \beta_j^2}$ and $\epsilon = (-1)^{\frac{p^n-1}{2}}$, so that \[ \det(M(\varphi)) = \epsilon \displaystyle\sum_{i \in \ZZ / p^n \ZZ} c_i \varphi_i^{p^n}. \]

Now observe that, taking $\varphi = e_l^*$ for some $l \in \ZZ/ p^n \ZZ$, the matrix $M(e_l^*)$ is monomial~: its term of index $(i,j)$ is $\delta_{l,i+j} \alpha_{ij}$.  Its determinant is $\det(M(e_l^*))= \epsilon \displaystyle\prod_{i+j=l} \alpha_{ij}$. We thus see that for all $l \in \ZZ / p^n \ZZ$ we have \[ c_l = \displaystyle\prod_{i+j = l} \alpha_{ij} .\]

Finally we can conclude that, for all $\varphi \in A^*$, we have 

\[ \det(M(\varphi)) = \epsilon \displaystyle\sum_{l \in \ZZ / p^n \ZZ} (\displaystyle\prod_{i+j = l} \alpha_{ij}) \varphi_l^{p^n}. \]

The $A_0$-algebra $A$ is Gorenstein if and only if there exists $\varphi \in A^*$ such that $\det(M(\varphi))$ is invertible. Since $A_0$ is a local ring, such is the case if and only if there exists $l \in \ZZ/ p^n \ZZ $ such that $\alpha_{ij}$ is invertible whenever $i+j = l$, in which case we can take $e_l^*$ as a generator of $A^*$.

\end{demo}

\subsection{Application to a Riemann-Hurwitz-type formula.}

We can now turn to the main main result of this section, which will relate the ramification divisor of a $\D(M)$-covering to the dualizing sheaf of its quotient morphism. We will use a formula known in height one (under the extra assumption that the base is an algebraically closed filed of characteristic $p>0$) and extend it to arbitrary height via our formalism.

Let us fix an algebraically closed field $k$ of characteristic $p>0$ and set $S=\Spec(k)$.

Let us first quote the following recent result (see also \cite{RS}) :

\begin{theo} \cite[Th 8.1]{Tziolas}
\label{Tzi}
Let $Y$ be an integral $S$-scheme. Suppose $Y$ has a dualizing sheaf $\omega_Y$, satisfies Serre's $S_2$ condition and has at worst normal crossing singularities in codimension $1$. Suppose $Y$ admits a $\mu_{p,S}$-action. Let $\cI_{\fix}$ be the ideal sheaf defining the scheme of fixed points and $f : Y \too X$ be the quotient. Then $X$ has a dualizing sheaf $\omega_X$ and \[ \omega_{Y} = ( f^* \omega_X \otimes \cI_{\fix}^{[1-p]})^{[1]}, \]

where, for every $\cO_Y$-module $\cF$, we denoted by $\cF^{[n]} = (\cF^{\otimes n})^{**}$ its $n$-th reflexive power.  
\end{theo}

Now suppose that a morphism $f : Y \too X$ is the quotient morphism of a $\mu_{p^n,S}$-covering of the scheme $X$. Under the extra assumption that $f$ is Gorenstein, which by the previous section can be checked on the cocycle giving the action on $Y$, we will prove an equality between its dualizing sheaf and the structure sheaf of the ramification divisor.

\begin{theo}

Let $X$ be a noetherian $S$-scheme and $(f : Y \too X, \mu_{p^n,S})$ be a $\mu_{p^n,S}$-covering of $X$. Denote by $\R_\cG$ the ramification divisor of the covering, as defined in \ref{newdiv}, ie the divisor associated with the action groupoid $\cG:=\mu_{p^n,S} \times Y \toto Y$. Suppose that $f$ is Gorenstein and denote by $\omega_f$ its dualizing sheaf. 

We then have \[ \omega_f = \cO_Y(\R_\cG). \]

\label{newRH}

\end{theo}

\begin{demo}

We will show the result by induction on $n$, using Tziolas's result. 

For $n=1$ this is a direct consequence of \ref{Tzi}. Indeed, since by assumption $f$ is Gorenstein, $\omega_f$ is invertible so in particular reflexive. Since $(\cI_{\fix}^{\otimes p-1})^*$ is the dual of a finite type module, it is also reflexive. Now by the previous result we have $(\omega_f)^* = \cI_{\fix}^{[p-1]}$ and hence $\omega_f = (\cI_{\fix}^{\otimes p-1})^*$. Furthermore since $Y$ is noetherian and $\cI_{\fix}^{\otimes p-1}$ is coherent by \cite[III, prop. 6.8]{Hartshorne} for all point $y \in Y$ we have $(\cI_{{\fix},y}^{\otimes p-1})^*= (\cI_{\fix}^{p-1})_y^*$. But by \ref{comparaison} for all point $y \in Y$ of codimension $1$ we have $\cO_Y(-\R_\cG)_y = \cI_{\fix,y}^{\otimes p-1}$. Thus dualizing we get $\cO_Y(\R_\cG)_y = (\cI_{{\fix},y}^{\otimes p-1})^*= (\cI_{\fix}^{p-1})_y^* = \omega_{f,y}$. Hence the sheaves $\cO_Y(\R_\cG)$ and $\omega_f$ are invertible sheaves equal in all codimension $1$ points so they are equal.    

Let $n>1$ be an integer and suppose the result is proved for all $\mu_{p^m,S}$-coverings with $m<n$. Let $(f : Z \too X, \mu_{p^n,S})$ be a covering of $X$. Consider $\mu_{p^{n-1},S}$ as a subgroup of $\mu_{p^n,S}$ via the obvious closed immersion and the induced action on $Z$. Set $Y = Z/\mu_{p^{n-1},S}$, which is normal since $Z$ is. We have a $\mu_{p^{n-1},S}$-covering $(g : Z \too Y, \mu_{p^{n-1},S})$ of $Y$ a commutative diagram  \[ \xymatrix{ Z \ar[d]_f \ar[r]^g & Y \ar[ld]^h \\ X} \] from which we see that $(h : Y \too X, \mu_{p,S})$ is a covering of $X$. If $\cG$, $\cH$, $\cQ$ are the associated action groupoids, by \ref{devissage} we have $\R_\cG = \R_\cH + g^*\R_\cQ$. By induction hypothesis we have $\cO_Y(\R_\cQ) = \omega_h$ and $\cO_Z(\R_\cH) = \omega_g$. But by \cite[6.4, Lemma 4.26]{Liu} we have $\omega_f = \omega_g \otimes_{\cO_Z} g^* \omega_h$. Thus \[ \omega_f = \cO_Z(\R_\cH) \otimes_{\cO_Z} g^*\cO_Y(\R_\cQ)= \cO_Z(\R_\cH + g^* \R_\cQ) = \cO_Z(\R_\cG) \] which shows that the formula holds for $\mu_{p^n,S}$.      
\end{demo}

The above formula extends immediately to the case of coverings given by actions of an arbitrary finite diagonalisable group scheme. Indeed, if $G = \D(M)$ is a finite  infinitesimal group scheme, the decomposition of $M$ into invariant factors yields a decomposition of $G$ into a product $G = \displaystyle\prod\limits_{i=1}^r \mu_{p^{n_i},S} \times G_{\textup{\'et}}$, where $G_{\textup{\'et}}$ is a finite \'etale $S$-group scheme. If $(f : Y \too X, G)$ is a $G$-covering of $X$ we can decompose the action of $G$ accordingly and factor $f$ into \[Y=Y_0 \stackrel{f_1}{\too} Y_1 \stackrel{f_2}{\too} \dots \stackrel{f_{r-1}}{\too} Y_r \stackrel{g}{\too} Y_{r+1}=X \] where $f_i$ is the quotient by $\mu_{p^{n_i},S}$ and $G$ is the quotient by $G_{\textup{\'et}}$. Successive applications of the above result \ref{newRH}, along with the classical formula (\ref{RH}) yields the following theorem :      

\begin{theo}

\label{newnRH}
Let $X$ be a scheme defined over an algebraically closed field $k$ and $G$ be a finite diagonalizable $k$-group scheme. Let $(f : Y \too X, G)$ be a covering of $X$ given by an action of $G$ on $Y$. Let $\R_\cG$ be the ramification divisor of this covering, associated with the action groupoid $\cG$.  If $f$ is Gorenstein then we have \[\omega_f = \cO_Y(\R_\cG), \]

where $\omega_f$ is the dualizing sheaf of the morphism $f$.

\end{theo}

As in the classical case, when $Y$ is a smooth projective curve over $k$ we can take the degrees in the above formula to relate the genuses of $Y$,$X$ and the degree of $\R_\cG$, as follows.

\begin{coro}

With the notations of the above theorem, suppose furthermore that $Y$ (and hence $X$) is a smooth projective curve over $k$. Let $g(Y)$, respectively $g(X)$, denote the genuses of the curve $Y$, respectively $X$. We have the formula \[ 2g(Y) - 2 =    \vert G \vert (2g(X)-2) + \deg(\R_\cG). \]

\end{coro}

\begin{demo}

By \cite[6.4, Lemma 4.26]{Liu} we have $\omega_{Y/k} \simeq f^*\omega_{X/k} \otimes_{\cO_Y} \omega_f$. Taking degrees we get $\deg(\omega_{Y/k}) = \deg(f) \deg(\omega_{X/k}) + \deg(\omega_f)$. By \ref{newnRH} we have $\deg(\omega_f) = \deg(\R_\cG)$.  Since the degree of the canonical divisor of a smooth projective curve $C$ is $2g(C)-2$, we get the announced formula, noting that $f$ is finite flat of degree $\vert G \vert$. 
\end{demo}

Note that this formula was proved by Emsalem in \cite[cor. 7.3]{emsalem} in the special case of torsors.

\medskip

\section*{Acknowledgements}

The author wishes to thank Matthieu Romagny for his constant support during all stages of this work. We would also like to thank Jo\~ao Pedro Dos Santos for his insightful remarks on a preliminary version of this article.

\bibliographystyle{amshalpha}
 \bibliography{biblio2}

 \end{document}